\documentclass{amsart}
\usepackage{mathrsfs}
\usepackage{amsfonts,enumerate,amsmath,amssymb}
\usepackage{color,latexsym,amsfonts,amssymb}

\newtheorem{thm}{Theorem}[section]

\newtheorem{lem}[thm]{Lemma}

\newtheorem{rem}[thm]{Remark}
\newtheorem{defn}[thm]{Definition}

\def\para#1{\vskip .4\baselineskip\noindent{\bf #1}}

\theoremstyle{definition}

\theoremstyle{remark}

\numberwithin{equation}{section}

\begin{document}
\title[Averaging for Mixed Fast-Slow Systems with FBm]{Averaging Principles for Mixed Fast-Slow Systems Driven by Fractional Brownian Motion}
\author{Bin Pei}
\address{School of Mathematics and Statistics, Northwestern Polytechnical University, Xi'an, 710072, China, Innovation Center, NPU-Chongqing, Chongqing, 401135, China, and Faculty of Mathematics, Kyushu University, Fukuoka, 819-0395, Japan}
\email{binpei@nwpu.edu.cn}

\author{Yuzuru Inahama}
\address{Faculty of Mathematics,
	Kyushu University, Fukuoka, 819-0395, Japan}
\email{inahama@math.kyushu-u.ac.jp}

\author{Yong Xu}
\address{Corresponding author, School of Mathematics and Statistics, Northwestern Polytechnical University, Xi'an, 710072, China}
\email{hsux3@nwpu.edu.cn}

\begin{abstract}
We focus on fast-slow systems involving both fractional Brownian motion (fBm) and standard Brownian motion (Bm). The integral with respect to Bm is the standard It\^{o} integral, and the integral with respect to fBm is a generalised Riemann-Stieltjes integral by means of fractional calculus. We establish an averaging principle in which the fast-varying diffusion process of the fast-slow systems acts as a \textquotedblleft noise\textquotedblright~to be averaged out in the limit. We show that the slow process has a limit in the mean square sense, which is characterized by the solution of stochastic differential equations driven by fBm whose coefficients are averaged with respect to the stationary measure of the fast-varying diffusion. An implication is that one can ignore the complex original systems and concentrate on the averaged systems instead. This averaging principle paves the way for reduction of computational complexity.
\end{abstract}
\maketitle
\section{Introduction}\label{se-1}
The real-valued fractional Brownian motion (fBm) with Hurst index $H\in (0,1)$ is a zero mean Gaussian process $\{B_{t}^{H}, t \geq 0\}$ with covariance  function
\begin{eqnarray}\label{RH}
\mathbb{E}[B_{t}^{H}B_{s}^{H}]=\frac{1}{2}(t^{2 H}+s^{2 H}-|t-s|^{2 H}).
\end{eqnarray}
From (\ref{RH}), we deduce that $
\mathbb{E}[|B_{t}^{H}-B_{s}^{H}|^{2}]=|t-s|^{2 H}$,  as a consequence, the trajectories of $B^H$ are almost surely locally $H'$-H\"{o}lder continuous for all $H' \in(0,H)$. This process was introduced by Kolmogorov \cite{Kolmogorov1940} and later studied by Mandelbrot and Van Ness \cite{Mandelbrot1968}. Its self-similar and long-range dependence ($H>\frac{1}{2}$) properties make this process a useful driving noise in models arising in physics, finance and other fields \cite{Biagini2008stochastic,decreusefond1998fractional,mishura2008stochastic}. Since $B^H$ is not a semimartingale if $H\neq \frac{1}{2}$, we cannot use the classical It\^{o} theory to construct a stochastic calculus with respect to the fBm. The $m$-dimensional fBm with same Hurst index $H$ is just a collection of $m$-independent one-dimensional fBm's in (\ref{RH}).

Over the last years some new techniques have been developed in order to define stochastic integrals with respect to fBm, see e.g. \cite{Biagini2008stochastic,Inahama2019,mishura2008stochastic,Nualart2006}. Lyons \cite{Lyons1994} solved the equations driven by a fBm with Hurst parameter $H>\frac{1}{2}$ by a pathwise approach using the $p$-variation norm. Nualart and R\u{a}\c{s}canu \cite{Rascanu2002} studied the differential equations driven by fBm using the tools of fractional calculus in the sense of Z\"{a}hle \cite{Zahle1998}. Kubilius \cite{Kubilius2002} studied one dimensional stochastic differential equations (SDEs) driven by both fBm and standard Brownian motion (Bm), with the noise term independent of the time and with no drift term.  Guerra and Nualart \cite{Guerra2008} proved an existence and uniqueness theorem for solutions of multidimensional, time dependent, SDEs driven by fBm with Hurst parameter $H>\frac{1}{2}$ and standard Bm.

Let $(\Omega, \mathscr{F}, \{\mathscr{F}_t\}_{t\geq 0},\mathbb{P})$ be a stochastic basis satisfying the usual conditions.  Take an aribitrary $H\in (\frac{1}{2},1)$ and fix it throughout this paper.  Let $B^H=\{B^H_t, t\in [0,T]\} $ and $W =\{W_t, t\in [0,T]\} $  be independent $m$-dimensional fBm adapted to $\{\mathscr{F}_t\}$ with Hurst parameter $H$ and $r$-dimensional $\{\mathscr{F}_t\}$-Bm, respectively.

This paper will consider the following mixed fast-slow systems driven by fBm:
\begin{eqnarray}\label{orinteg1}
\left\{\begin{array}{l}{d X_{t}^{\varepsilon}=b_1(t, X_{t}^{\varepsilon}, Y_{t}^{\varepsilon}) d t+\sigma_1(t,X^\varepsilon_t) d B^{H}_t,} \quad \quad X_{0}^{\varepsilon}=x_0, \\ {d Y_{t}^{\varepsilon}=\frac{1}{\varepsilon} b_2( X_{t}^{\varepsilon}, Y_{t}^{\varepsilon}) d t+\frac{1}{\sqrt{\varepsilon}} \sigma_{2}(X_{t}^{\varepsilon}, Y_{t}^{\varepsilon}) d W_{t},}\quad Y_{0}^{\varepsilon}=y_0,\end{array}\right.
\end{eqnarray}
where the parameter $0<\varepsilon \ll 1$ represents the ratio between the natural time scale of the variables $X_{t}^{\varepsilon}$ and $Y_{t}^{\varepsilon}$ and $x_0 \in \mathbb{R}^{d_1}$ and $y_0 \in \mathbb{R}^{d_2}$ are arbitrary and non-random but fixed and the coefficients are measurable functions  $b_1^{l_1}: [0,T] \times  \mathbb{R}^{d_1}  \times \mathbb{R}^{d_2}\rightarrow \mathbb{R}, \sigma_{1}^{l_1, j}: [0,T] \times \mathbb{R}^{d_1} \rightarrow \mathbb{R}, 1 \leq l_1 \leq d_1, 1 \leq j \leq m$ and $b_2^{i},\sigma_{2}^{i, l_2}: \mathbb{R}^{d_1}  \times \mathbb{R}^{d_2}\rightarrow \mathbb{R}, 1 \leq i \leq d_2,  1 \leq l_2 \leq r$.  The integral $\int  \cdot~d W $ should be interpreted as an It\^{o} stochastic integral and the integral $\int \cdot ~d B^{H}$ as a generalised Riemann-Stieltjes integral in the sense of  Z\"{a}hle \cite{Zahle1998,Rascanu2002,Guerra2008}.

We will make use of the following assumptions (H1) and (H2) to obtain the existence and uniqueness result to  (\ref{orinteg1}).
\begin{description}
	\item[(H1)] The function
	$\sigma_1(t,x)$ is continuous and continuously
	differentiable in the variable $x$ and H\"older continuous in $t$ and has linear growth in the variable $x$, uniformly in $t$. Precisely, there exist constants
	$L_i,i=1,2,3,4$, and some constants $0 < \beta,\gamma \leq1$, such that
	\begin{eqnarray*}
		|\nabla_x \sigma_1(t,x_1)| &\leq& L_{1}\cr
		|\nabla_x \sigma_1(t,x_1)-\nabla_x \sigma_1(t,x_2)| &\leq & L_{2}|x_1-x_2|^{\gamma},\cr
		|\nabla_x \sigma_1(t,x_1)-\nabla_x \sigma_1(s,x_1)|+|\sigma_1(t,x_1)-\sigma_1(s,x_1)|
		&\leq& L_{3}|t-s|^{\beta},\cr
		|\sigma_1(t,x_1)|&\leq& L_{4}(1+|x_1|)
	\end{eqnarray*}
	for any $x_1,x_2\in \mathbb{R}^{d_1}$ and $t,s \in [0,T].$ Here, $\nabla_x$  is the standard gradient with respect to the $x$-variable.
	
	\item [(H2)] The function $b_1(t,x,y)$ is Lipschitz continuous in the variables $t, x,y$ and has linear growth in the variables $x,y$, uniformly in $t$ and the functions $b_2(x,y),$ $ \sigma_{2}(x,y)$ are Lipschitz continuous in the variables $x,y$ and have linear growth in the same variables. Precisely, there exist constants $L_i,i=5,6,7$, such that
	\begin{eqnarray*}
		|b_1(t,x_1,y_1)-b_1(s,x_2,y_2)|&\leq& L_5(|x_1-x_2|\cr
		&&+|y_1-y_2|+|t-s|),\cr
		|b_2(x_1,y_1)-b_2(x_2,y_2)|+|\sigma_{2}(x_1,y_1)-\sigma_{2}(x_2,y_2)|
		&\leq& L_6(|x_1-x_2|+|y_1-y_2|),\cr
		|b_1(t,x_1,y_1)|+|b_2(x_1,y_1)|+|\sigma_{2}(x_1,y_1)| &\leq& L_7(1+|x_1|+|y_1|)
	\end{eqnarray*}
	for any $x_1, x_2 \in \mathbb{R}^{d_1}, y_1, y_2 \in \mathbb{R}^{d_2}$ and $t,s \in [0,T].$
\end{description}

Now, we define the averaged equation:
\begin{eqnarray}\label{xave}
d\bar{X}_{t}=\bar{b}_1(t,\bar{X}_{t}) d t+\sigma_1(t,\bar{X}_{t}) d B^{H}_t, \ \ \ \bar{X}_{0}=x_0,
\end{eqnarray}
where $$\bar{b}_1(t,x)=\int_{\mathbb{R}^{d_2}} b_1(t,x, y) \mu^{x}(d y), \ \ \ x \in \mathbb{R}^{d_1},$$ and $\mu^{x} $ is a unique invariant probability measure which will be given in Appendix A (see also \cite[Proposition 3.8]{Liu2019} for example) with respect to the following frozen equation (\ref{frozon}) under conditions (H2) and (H4) below.
\begin{eqnarray}\label{frozon}
dY^{x}_t=b_2(x,Y^{x}_t)dt+\sigma_{2}(x,Y^{x}_t)d{W}_t, \quad Y^{x}_0=y.
\end{eqnarray}

To establish the averaging principle of  (\ref{orinteg1}), we set the following assumptions:
\begin{description}
	\item[(H3)] Assume further that $\sup_{t \in [0,T],x\in \mathbb{R}^{d_1},y \in \mathbb{R}^{d_2}}|b_1(t,x,y)|<\infty$ holds.
	\item [(H4)] There exist $C>0, \beta_i>0,i=1,2$, such that
\begin{align*}
		2 \langle y_1-y_2, b_2(x,y_1)-b_2(x,y_2) \rangle +&|\sigma_{2}(x,y_1)-\sigma_{2}(x,y_2)|^2 \\
		&\leq -\beta_1 |y_1-y_2|^2,\\
2 \langle y_1, b_2(x,y_1) \rangle +|\sigma_{2}(x,y_1)|^2&\leq -\beta_2 |y_1|^2+C|x|^2+C,
\end{align*}
	hold for any $x\in \mathbb{R}^{d_1}$, $y_1,y_2\in \mathbb{R}^{d_2}$.
\end{description}

We follow the approach by \cite{Guerra2008,Rascanu2002} and introduce some necessary spaces and norms.
Taking a parameter $ 1-H<\alpha<\frac{1}{2}$, denote by $W_{0}^{\alpha,\infty}$ the space of measurable functions $f:[0,T] \rightarrow\mathbb{R}^d$  such
that
$$\|f\|_{\alpha, \infty} :=\sup _{t \in[0, T]}\|f(t)\|_{\alpha}<\infty,$$
where $$\|f(t)\|_{\alpha}=|f(t)|+\int_{0}^{t} \frac{|f(t)-f(s)|}{(t-s)^{\alpha+1}} d s.$$

Note that $C$ and $C_{{\rm x}}$ denote some positive constants which may change from line to line throughout this paper, where ${\rm x}$ is one or more than one parameter and $C_{{\rm x}}$ is used to emphasize that the constant depends on the corresponding parameter, for example, $C_{\alpha,\beta,\gamma, T, R, |x_0|,|y_0|}$ depends on $\alpha,\beta,\gamma, T, R, |x_0|$ and $|y_0|$.

Then, we formulate our main result of averaging principle in the mean square sense.
\begin{thm}\label{thm1}
	Suppose that {\rm (H1)-(H4)} hold and let $\beta$ and $\gamma$ be as in {\rm (H1)}. Let $1-H<\alpha <\min \{\frac{1}{2}, \beta,\frac{\gamma}{2}\}$, then we obtain
	\begin{eqnarray*}
		\lim\limits_{\varepsilon \rightarrow 0} \mathbb{E}\big[\|{X}^{\varepsilon}-\bar{X}\|^2_{\alpha,\infty}\big] =0.
	\end{eqnarray*}
\end{thm}

\begin{rem}{\rm
		From Theorem \ref{thm1}, we know that, for different $\beta$ and $\gamma$, ${X}^{\varepsilon}$ converges to $\bar{X}$ in the sense of mean square with different Hurst exponents, i.e.,
		\begin{itemize}
			\item $H\in(\frac{1}{2},1)$, \ \ \ \ \ \ if $\gamma=1,\beta\geq \frac{1}{2}$;
			\item $H\in (1-\frac{\gamma}{2},1)$,  \  if $\gamma<1,\beta \geq \frac{\gamma}{2}$;
			\item $H\in (1-\beta,1)$, \ if $\beta < \min\{\frac{\gamma}{2},\frac{1}{2}\}$.
		\end{itemize}
		
		In particular, if $\sigma_1$ is regular enough, our main theorem can be applied for any $H\in (\frac{1}{2},1)$.}
\end{rem}		

\begin{rem}{\rm
		In order to obtain the strong convergence, it is known that the diffusion coefficient $\sigma_1$ in (\ref{orinteg1}) should not depend on the fast variable $Y^\varepsilon$ (see e.g. \cite{Givon2007}). }
\end{rem}	

The study on averaging principles for stochastic systems can be traced back to the work of
Khasminskii \cite{khas1966limit}, see also the recent effort in \cite{duan2014,Duan2015,fu2015strong,KY1,KY3,Liu2019,peiconv2020,peisto2020,peisto2017,xu2014,xu2015approximation,Xu2015Stochastic,Xu2017Stochastic,wu2020} and references therein. A central theme is:
The fast varying process can be treated as a noise and has an invariant measure. Utilizing this invariant measure, one can carry out asymptotic analysis so
that the slow process converges
to a limit that is an average with respect to the stationary measure of the fast-varying process.
Freidlin
and Wentzell \cite{freidlin2012random} provided an illuminating overview and discussion on the averaging principle. Givon \cite{Givon2007} studied the two-time-scale jump-diffusion stochastic differential systems and obtained the strong convergence rate of the slow components to the effective dynamics. Thompson, Kuske and Monahan \cite{thompson2015stochastic} studied nonlinear fast-slow stochastic dynamical systems in which the fast variables are driven by additive $\alpha$-stable noise perturbations and the slow variables depending linearly on the fast variables. Xu and Miao \cite{xumiao2015} studied $L^p$-strong convergence of an averaging principle for two-time-scales jump-diffusion SDEs. Xu, Pei and Guo \cite{Xu2015Stochastic} investigated the stochastic averaging of slow-fast dynamical systems driven by fBm with the Hurst parameter $H$ in the interval $(\frac{1}{2},1)$. Hairer and Li \cite{Hairer2019} considered slow-fast systems where the slow system is driven by fBm and proved the convergence to the averaged solution took place in probability which strongly relies on stochastic sewing lemma.

Nevertheless, the aforementioned paper cannot answer the question that if disturbances involve both standard Bm and long-range dependence modeled by fBm $H\in (\frac{1}{2},1)$ in the mean square sense. In this paper, we aim to address this issue.
We answer affirmatively that an averaging principle still holds for fast-slow systems involving both standard Bm and fBm. The main difficulties here are how to deal with fBm, standard Bm. In order to overcome these difficulties, our approach is completely different from Xu's previous work \cite{Xu2015Stochastic} in the sense that we combine the pathwise approach with the It\^{o} stochastic calculus to handle both types of integrals and use stopping time techniques to establish averaging principle for multidimensional, time dependent, SDEs driven by fBm with fast-varying diffusion process.

The rest of the paper is organized as follows. Section \ref{se-2} presents some necessary notations and assumptions. The existence and uniqueness theorems for (\ref{orinteg1}) and (\ref{xave}) were proved in Section \ref{sec-3}. Section \ref{se-4} presents fast-slow systems driven by fBm with fast-varying diffusion process. Stochastic averaging principles for such SDEs are then established. Some technical complements are included in the appendix, which provides the arguments of the ergodicity for the fast component in which the slow component is kept frozen.

\section{Preliminaries}\label{se-2}
We recall some basic facts on generalised Riemann-Stieltjes integrals. For more details, we refer to the paper \cite{Rascanu2002,Guerra2008} and a monograph \cite{mishura2008stochastic}. Let $f \in L^1(a,b)$ and $\alpha >0$. The fractional left and right Riemann-Liouville integrals of order $\alpha$ are defined for almost all $x\in(a,b)$ by
$$
I_{a+}^{\alpha} f(x)=\frac{1}{\Gamma(\alpha)} \int_{a}^{x} \frac{f(y)}{(x-y)^{1-\alpha}} d y,
$$
and
$$
I_{b-}^{\alpha} f(x)=\frac{(-1)^{-\alpha}}{\Gamma(\alpha)} \int_{x}^{b} \frac{f(y)}{(y-x)^{1-\alpha}} d y,
$$
respectively, where $(-1)^\alpha=e^{-i\pi \alpha}$ and $\Gamma(\alpha)=\int_{0}^{\infty}r^{\alpha-1}e^{-r}dr$ is the Euler Gamma function. Let $I^{\alpha}_{a+}(L^p)$ (resp. $I^{\alpha}_{b-}(L^p)$) be the image of $L^{p}(a,b)$ by the operator $I^{\alpha}_{a+}$ (resp. $I^{\alpha}_{b-}$). If $f\in I^{\alpha}_{a+}(L^p)$ (resp. $f\in I^{\alpha}_{b-}(L^p)$) and $0<\alpha <1$, then
the Weyl derivatives of $f$ are defined by formulas
$$
D_{a+}^{\alpha} f(x):=\frac{1}{\Gamma(1-\alpha)}\bigg(\frac{f(x)}{(x-a)^{\alpha}}+\alpha \int_{a}^{x} \frac{f(x)- f(y)}{(x-y)^{\alpha+1}} d y\bigg)\mathbf{1}_{(a,b)}(x),
$$
and $$
D_{b-}^{\alpha} f(x):=\frac{(-1)^\alpha}{\Gamma(1-\alpha)}\bigg(\frac{f(x)}{(b-x)^{\alpha}}+\alpha \int_{x}^{b} \frac{f(x)-f(y)}{(y-x)^{\alpha+1}} d y\bigg)\mathbf{1}_{(a,b)}(x),
$$
and are defined for almost all $x\in(a,b)$ (the convergence of the integrals at the singularity $y=x$ holds pointwise for almost all $x\in(a,b)$ if $p=1$ and moreover in $L^p$-sense if $1<p<\infty$).

We have that:
\begin{itemize}
	\item If $\alpha<\frac{1}{p}$ and $q=\frac{p}{1-\alpha p}$, then $$I^{\alpha}_{a+}(L^p)=I^{\alpha}_{b-}(L^p) \subset L^{q} (a,b).$$
	\item If $\alpha>\frac{1}{p}$, then $$I^{\alpha}_{a+}(L^p) \cup I^{\alpha}_{b-}(L^p) \subset C^{\alpha-\frac{1}{p}} (a,b).$$
\end{itemize}
The fractional integrals and derivatives are related by the inversion formulas
\begin{eqnarray*}
	&&I_{a{+}}^{\alpha}(D_{a{+}}^{\alpha} f) =f, \quad \forall f \in I_{a{+}}^{\alpha}(L^{p}), \cr
	&&D_{a{+}}^{\alpha}(I_{a{+}}^{\alpha} f) =f, \quad \forall f \in L^{1}(a, b),
\end{eqnarray*}
and similar statements also hold for $I^{\alpha}_{b-}$ and $D^{\alpha}_{b-}$.

Let $f(a+) :=\lim _{\varepsilon \searrow 0} f(a+\varepsilon)$ and $g(b-) :=\lim _{\varepsilon \searrow 0} g(b-\varepsilon)$  (we are assuming that these limits exist and are finite) and define
\begin{eqnarray*} &&f_{a{+}}(x) :=(f(x)-f(a+)) \mathbf{1}_{(a, b)}(x), \cr &&g_{b{-}}(x) :=(g(x)-g(b-)) \mathbf{1}_{(a, b)}(x). \end{eqnarray*}
We recall from Z\"{a}hle \cite{Zahle1998}, the definition of generalized Riemann-Stieltjes
fractional integral with respect to irregular functions.
\begin{defn}{\rm
		(Generalized Riemann-Stieltjes Integral). Let $f$ and $g$ be functions such that the limits
		$f(a+), g(a+), g(b-)$ exist. Suppose that $f_{a+} \in I_{a+}^{\alpha}(L^{p})$ and $g_{b-} \in I_{b-}^{1-\alpha}(L^{q})$ for some $\alpha \in(0,1)$ and $p, q \in[1, \infty] $ such that
		$\frac{1}{p}+\frac{1}{q}\leq1$. In this case the generalised Riemann-Stieltjes integral
		$$\begin{aligned} \int_{a}^{b} f d g =(-1)^{\alpha} \int_{a}^{b} D_{a+}^{\alpha}f_{a+}(x) D_{b-}^{1-\alpha}g_{b-}(x)dx+f(a+)(g(b-)-g(a+)), \end{aligned}
		$$
		is well-defined. }
\end{defn}

For $\eta\in(0,1]$, let $C^\eta,$ be the space of $\eta$-H\"{o}lder continuous functions $f:[a,b]\rightarrow\mathbb{R}^d,$ equipped with the the norm $$
\|f\|_{\eta} :=\|f\|_{\infty}+\sup_{a \leq s < t \leq b}\frac{|f(t)-f(s)|}{(t-s)^{\eta}}<\infty,$$
where $\|f\|_{\infty}=\sup _{t \in[a, b]}|f(t)|.$
Given any $\epsilon$ such that $0<\epsilon<\alpha$, we have the following continuous inclusions $C^{\alpha+\epsilon}\subset W_{0}^{\alpha,\infty}\subset C^{\alpha-\epsilon}.$

\begin{rem}
	{\rm The above definition is simpler in the following cases.
		\begin{itemize}
			\item If $\alpha<\frac{1}{p}$, under the assumptions of the preceding definition, we have
			that $f\in I^{\alpha}_{a+}(L^p)$ and we can write
			\begin{eqnarray}\label{rsite}
			\int_{a}^{b} f {d} g=(-1)^{\alpha} \int_{a}^{b} D_{a{+}}^{\alpha} f(x) D_{b-}^{1-\alpha} g_{b-}(x) {d} x.
			\end{eqnarray}
			\item If $f\in C^{\eta_1}(a,b)$ and $g\in C^{\eta_2}(a,b)$ with $\eta_1+\eta_2>1$ then we can choose $p=q=\infty$ and  $1-\eta_2<\alpha<\eta_1$, the generalized Riemann-Stieltjes integral exists, it is given by (\ref{rsite}) and coincides with the Riemann-Stieltjes integral.
	\end{itemize}}
\end{rem}

Now, fix the parameter $\alpha$, such that $0<
\alpha<\frac{1}{2}$, denote by $W_{0}^{\alpha,1}$
the space of measurable functions  $f:[0,T] \rightarrow \mathbb{R}^d$
such that
$$
\|f\|_{\alpha, 1} :=\int_{0}^{T} \frac{|f(s)|}{s^{\alpha}} d s+\int_{0}^{T} \int_{0}^{s} \frac{|f(s)-f(y)|}{(s-y)^{\alpha+1}} d y d s<\infty.
$$

Denote by $W_{T}^{1-\alpha,\infty}$ the space of measurable functions  $g:[0,T]\rightarrow\mathbb{R}^m$ such that
$$\|g\|_{1-\alpha, \infty, T} :=\sup _{0<s<t<T}\bigg(\frac{|g(t)-g(s)|}{(t-s)^{1-\alpha}}+\int_{s}^{t} \frac{|g(y)-g(s)|}{(y-s)^{2-\alpha}} d y\bigg)<\infty.$$
It is also easy to prove that $C^{1-\alpha+\epsilon} \subset W_{T}^{1-\alpha,\infty} \subset C^{1-\alpha}$. For $g\in W_{T}^{1-\alpha,\infty}$, we have that
\begin{eqnarray*}
	\Lambda_{\alpha}(g) &:=&\frac{1}{\Gamma(1-\alpha)} \sup _{0 < s<t <T}|(D_{t-}^{1-\alpha} g_{t-})(s)| \cr
	& \leq & \frac{1}{\Gamma(1-\alpha) \Gamma(\alpha)}\|g\|_{1-\alpha, \infty, T}<\infty.
\end{eqnarray*}
Moreover, if $f\in W_{0}^{\alpha,1}$ and $g\in W_{T}^{1-\alpha,\infty}$  then $\int_{0}^{t}fdg$ exists for all $t\in[0,T]$ and
$$\bigg|\int_{0}^{t} f d g\bigg| \leq \Lambda_{\alpha}(g)\|f\|_{\alpha, 1},$$
holds.

\begin{rem}\label{fbmito}
	{\rm The trajectories of $B^H$ are almost surely locally $H'$-H\"{o}lder
		continuous for all  $H' \in (0,H)$. Then, for all $1-H<\alpha <\frac{1}{2}$, the trajectories of $B^H$ belong to the space $W_{T}^{1-\alpha,\infty}$. As a consequence, the generalised  Riemann-Stieltjes integrals $$\int_{0}^{T} v_s dB^H_s$$ exists if $\{v_t,t\in[0,T]\}$ is a stochastic process whose trajectories belong to the space $W_{0}^{\alpha,1}$. And we have
		\begin{eqnarray}\label{fbm}
		\bigg|\int_{0}^{t} v_s d B_s^H\bigg| \leq \Lambda_{\alpha}(B^H)\|v\|_{\alpha, 1},
		\end{eqnarray}
		where $\Lambda_{\alpha}(B^H):=\frac{1}{\Gamma(1-\alpha) \Gamma(\alpha)}\|B^H\|_{1-\alpha, \infty, T}$ has moments of all order, see Lemma 7.5 in Nualart and R\u{a}\c{s}canu \cite{Rascanu2002}. Furthermore, by the classical Fernique's theorem, for any $0<\vartheta <2$, we have
		\begin{eqnarray}\label{expfbm}
		\mathbb{E}[\exp(\Lambda_{\alpha}^\vartheta(B^H))]<\infty.
		\end{eqnarray}}
\end{rem}

\section{Existence and Uniqueness for the Fast-slow Systems}\label{sec-3}
According to Theorem 2.2 in \cite{Guerra2008}, we obtain the existence and uniqueness result to  (\ref{orinteg1}).
\begin{lem}
	Suppose that {\rm (H1) and (H2)} hold and let $ 1-H<\alpha<\min\{\frac{1}{2},\beta,\frac{\gamma}{2}\}$. Then,  {\rm (\ref{orinteg1})} has a pathwise unique strong solution $\{(X_t^\varepsilon,Y^\varepsilon_t),t\geq 0\}$, i.e.,
	\begin{eqnarray}\label{mildx}
	\left\{\begin{array}{l}{X_{t}^{\varepsilon}=x_0+\int_{0}^{t}b_1(s, X_{s}^{\varepsilon}, Y_{s}^{\varepsilon}) d s+ \int_{0}^{t}\sigma_1(s,X^\varepsilon_s) d B^{H}_s,} \\ { Y_{t}^{\varepsilon}=y_0+\frac{1}{\varepsilon} \int_{0}^{t}b_2( X_{s}^{\varepsilon}, Y_{s}^{\varepsilon}) d s+\frac{1}{\sqrt{\varepsilon}}\int_{0}^{t} \sigma_{2}(X_{s}^{\varepsilon}, Y_{s}^{\varepsilon}) d W_{s}.}\end{array}\right.
	\end{eqnarray}
\end{lem}

\begin{lem}
	Suppose that {\rm (H1), (H2) and (H4)} hold and let $ 1-H<\alpha<\min\{\frac{1}{2},\beta,\frac{\gamma}{2}\}$. Then,  {\rm (\ref{xave})} has a pathwise unique strong solution $\{\bar{X}_t,t\geq 0\}$.
\end{lem}
\para{Proof:} For any $x_1, x_2, x \in \mathbb{R}^{d_1}$ and any initial value $y\in \mathbb{R}^{d_2}$, by Lemma \ref{invariant1} and Lemma \ref{x1-x2} in Appendix A, we have
\begin{eqnarray}
|\bar{b}_1(t,x_1)-\bar{b}_1(t,x_2)|
&\leq&\bigg|\int_{\mathbb{R}^{d_2}} b_1(t,x_1, z)\mu^{x_1}(d z)-\mathbb{E}[b_1(t,x_1,Y^{x_1,y}_s)]\bigg|\cr
&&+\bigg|\int_{\mathbb{R}^{d_2}} b_1(t,x_2,z) \mu^{x_2}(d z)-\mathbb{E}[b_1(t,x_2,Y^{x_2,y}_s)]\bigg|\cr
&&+\big|\mathbb{E}[b_1(t,x_1,Y^{x_1,y}_s)-b_1(t,x_2,Y^{x_2,y}_s)]\big|\cr
&\leq&Ce^{-\beta_1 s}(1+|x_1|+|x_2|+|y|)+C|x_1-x_2|,
\end{eqnarray}
and
\begin{eqnarray}
|\bar{b}_1(t_1,x)-\bar{b}_1(t_2,x)|
&\leq&\bigg|\int_{\mathbb{R}^{d_2}} b_1(t_1,x, z)\mu^{x}(d z)-\mathbb{E}[b_1(t_1,x,Y^{x,y}_s)]\bigg|\cr
&&+\bigg|\int_{\mathbb{R}^{d_2}} b_1(t_2,x,z) \mu^{x}(d z)-\mathbb{E}[b_1(t_2,x,Y^{x,y}_s)]\bigg|\cr
&&+\big|\mathbb{E}[b_1(t_1,x,Y^{x,y}_s)-b_1(t_2, x,Y^{x,y}_s)]\big|\cr
&\leq&Ce^{-\beta_1 s}(1+|x|+|y|)+C|t_1-t_2|.
\end{eqnarray}
Let $s\rightarrow \infty$, then we obtain that $\bar{b}_1$ is Lipschitz continuous in $x$ and $t$, and
\begin{eqnarray}
|\bar{b}_1(t,x)|\leq \int_{\mathbb{R}^{d_2}}|b_1(t,x, z)|\mu^{x}(d z)\leq C(1+|x|).
\end{eqnarray}

So, $\bar{b}_1$ satisfies the growth condition. Thus, according to Theorem 2.2 in \cite{Guerra2008},  (\ref{xave}) has a unique strong solution. \qed
\section{Proof of Main Result}\label{se-4}
This section is devoted to proving Theorem \ref{thm1}. The proof consists of the following steps.

Firstly, we give some a priori estimates for the solution $(X^{\varepsilon},Y^{\varepsilon})$ to  (\ref{orinteg1}).

Secondly, following the discretization techniques inspired by Khasminskii in \cite{khas1966limit}, we introduce an auxiliary process $(\hat{X}^{\varepsilon},\hat{Y}^{\varepsilon})$ and divide $[0,T]$ into intervals depending of size $\delta<1$, where $\delta$ is a fixed positive number depending on $\varepsilon$ which will be chosen later. Then, we construct $\hat{Y}^{\varepsilon}$ with initial value $\hat{Y}_{0}^{\varepsilon}=y_0, $ and for $t\in [k\delta, \min\{(k+1)\delta,T\}]$,
$$
\hat{Y}_{t}^{\varepsilon}= \hat{Y}_{k\delta}^{\varepsilon}+\frac{1}{\varepsilon} \int_{k\delta}^{t}b_2(X_{k \delta}^{\varepsilon}, \hat{Y}_{s}^{\varepsilon}) d s+\frac{1}{\sqrt{\varepsilon}} \int_{k\delta}^{t} \sigma_{2}(X_{k \delta}^{\varepsilon}, \hat{Y}_{s}^{\varepsilon}) d W_{s},
$$
i.e.
$$
\hat{Y}_{t}^{\varepsilon}=y_0+\frac{1}{\varepsilon} \int_{0}^{t}b_2(X_{s(\delta)}^{\varepsilon}, \hat{Y}_{s}^{\varepsilon}) d s+\frac{1}{\sqrt{\varepsilon}}\int_{0}^{t} \sigma_{2}(X_{s (\delta)}^{\varepsilon}, \hat{Y}_{s}^{\varepsilon}) d W_{s},
$$
where $s(\delta)=\lfloor \frac{s}{\delta} \rfloor \delta$  is the nearest breakpoint preceding $s$. Also, we define the process $\hat{X}^{\varepsilon}$ with initial value $\hat{X}_{0}^{\varepsilon}=x_0,$ by
\begin{eqnarray}\label{xhat}
\hat{X}_{t}^{\varepsilon}=x_0+\int_{0}^{t} b_1(s(\delta), X_{s(\delta)}^{\varepsilon}, \hat{Y}_{s}^{\varepsilon}) d s+\int_{0}^{t} \sigma_1(s,X_{s}^{\varepsilon}) d B^H_{s},
\end{eqnarray}
and then, we can derive uniform bounds $\|X_{t}^{\varepsilon}-\hat{X}_{t}^{\varepsilon}\|_{\alpha}$.

Thirdly, based on the ergodic property of the frozen equation, we obtain appropriate control of $\|\hat{X}_{t}^{\varepsilon}-\bar{X}_{t}\|_{\alpha}$.

Finally, we can estimate $\|{X}_{t}^{\varepsilon}-\bar{X}_{t}\|_{\alpha}$.

\para{Step 1: A priori estimates for the solution $(X^{\varepsilon},Y^{\varepsilon})$.}
We use techniques similar to those used in \cite[Theorem 4.2]{Shevchenko2014mixed} to give a priori estimate for the solution $X^{\varepsilon}$.	
\begin{lem}\label{lemb}
	Suppose that {\rm (H1)-(H3)} hold. Then, for $t\in[0,T], p\geq 1$, we have
	\begin{eqnarray*}
		\mathbb{E}[\|X^\varepsilon\|_{\alpha,\infty}^{p}] \leq C_{\alpha,T,|x_0|}.
	\end{eqnarray*}
\end{lem}
\para{Proof:} For shortness, denote, $\Lambda:=\Lambda_{\alpha}(B^{H})\vee 1,$ and for any $\lambda \geq 1, $ let
$$\|f\|_{\lambda,t} :=\sup _{0 \leq s \leq t}e^{-\lambda s}|f(s)|,$$ and
\begin{eqnarray*}
	\|f\|_{1,\lambda,t} :=\sup_{0 \leq s \leq t}e^{-\lambda s}\int_{0}^{s} \frac{|f(s)-f(r)|}{(s-r)^{\alpha+1}}dr.
\end{eqnarray*}

By (H3) and (\ref{fbm}), we start by estimating $\|X^\varepsilon\|_{\lambda,t}$:
\begin{eqnarray}\label{in1}
\|X^\varepsilon\|_{\lambda,t}&=&\sup _{0 \leq s \leq t}e^{-\lambda s}\bigg|x_0+\int_{0}^{s}b_1(r, X_{r}^{\varepsilon}, Y_{r}^{\varepsilon}) d r+ \int_{0}^{s}\sigma_1(r,X^\varepsilon_r) d B^{H}_r\bigg|\cr
&\leq& C_{T,|x_0|} \Lambda \Big(1+\sup _{0 \leq s \leq t}\int_{0}^{s} e^{-\lambda (s-r)}(r^{-\alpha} \|X^\varepsilon\|_{\lambda,t} +\|X^\varepsilon\|_{1,\lambda,t})dr\Big)\cr
&\leq& K \Lambda(1+\lambda^{\alpha-1}\|X^\varepsilon\|_{\lambda, t}+\lambda^{-1}\|X^\varepsilon\|_{1, \lambda, t}),
\end{eqnarray}
with some constant $K$ (which is dependent on $|x_0|$ and can be assumed to be greater than 1 without loss of generality) and here, we have used the estimate
\begin{eqnarray}\label{inq1}
\int_{0}^{t} e^{-\lambda(t-r)} r^{-\alpha} d r
&=& \frac{1}{\lambda} \int_{0}^{\lambda t} e^{-y} \lambda^{\alpha}(\lambda t-y)^{-\alpha} d y \cr
&\leq& \lambda^{\alpha-1} \sup _{z>0} \int_{0}^{z} e^{-y}(z-y)^{-\alpha} d y \cr
&\leq& C \lambda^{\alpha-1}.
\end{eqnarray}

Furthermore, we estimate $\|X^\varepsilon\|_{1, \lambda, t}$. To complete this step,
$$\mathcal{C}:=\int_{0}^t{(t-s)^{-\alpha-1}}{\bigg|\int_{s}^{t} f(r) d B^H_{r}\bigg|}ds,$$ needs to be estimated in advance, here, $f:[0,T] \rightarrow\mathbb{R}^d$ is a measurable function.
Using Fubini's theorem, it is easy to get
\begin{eqnarray*}
	\mathcal{C} &\leq& \Lambda_{\alpha}(B^H) \bigg(\int_{0}^t\int_{s}^{t}(t-s)^{-\alpha-1}\frac{|f(r)|}{(r-s)^{\alpha}}drds\cr
	&&+\int_{0}^t\int_{s}^{t} \int_{s}^{r} (t-s)^{-\alpha-1}\frac{|f(r)-f(q)|}{(r-q)^{1+\alpha}}dqdrds\bigg)\cr
	&\leq& \Lambda_{\alpha}(B^H) \bigg(\int_{0}^t\int_{0}^{r}(t-s)^{-\alpha-1}{(r-s)^{-\alpha}}ds{|f(r)|}dr\cr
	&&+\int_{0}^t\int_{0}^{r} \int_{0}^{q} (t-s)^{-\alpha-1}ds\frac{|f(r)-f(q)|}{(r-q)^{1+\alpha}}dqdr\bigg).
\end{eqnarray*}
Then, by the substitution $s=r-(t-r)y$, we have
\begin{eqnarray*}
	\int_{0}^{r}(t-s)^{-\alpha-1}{(r-s)^{-\alpha}}ds=(t-r)^{-2\alpha}\int_{0}^{r/(t-r)}(1+q)^{-\alpha-1}{q^{-\alpha}}dq,
\end{eqnarray*}
and on the other hand,
\begin{eqnarray*}
	\int_{0}^{q} (t-s)^{-\alpha-1}ds=\alpha^{-1}[(t-q)^{-\alpha}-t^{-\alpha}]\leq \alpha^{-1}(t-q)^{-\alpha}.
\end{eqnarray*}
This yields that
\begin{eqnarray}\label{fbmterm1}
\mathcal{C}
&\leq& \Lambda_{\alpha}(B^H) \bigg(c_\alpha \int_{0}^t (t-r)^{-2\alpha}{|f(r)|}dr\cr
&&+\int_{0}^t\int_{0}^{r} (t-q)^{-\alpha}\frac{|f(r)-f(q)|}{(r-q)^{1+\alpha}}dqdr\bigg),
\end{eqnarray}
where $c_\alpha=\int_{0}^{\infty}(1+q)^{-\alpha-1}{q^{-\alpha}}dq=B(2\alpha,1-\alpha)$, $B(2\alpha,1-\alpha)$ is the Beta function.

Thus, by (\ref{fbmterm1}), we have
\begin{eqnarray}\label{in2}
\|X^{\varepsilon}\|_{1,\lambda,t} &=&\sup_{0 \leq s \leq t}e^{-\lambda s}\int_{0}^{s} {(s-r)^{-\alpha-1}}{\bigg|\int_{r}^{s}b_1(q, X_{q}^{\varepsilon}, Y_{q}^{\varepsilon}) dq \bigg|} dr\cr
&&+\sup_{0 \leq s \leq t}e^{-\lambda s}\int_{0}^{s} {(s-r)^{-\alpha-1}}{\bigg|\int_{r}^{s}\sigma_1(q,X^\varepsilon_q) d B^{H}_q\bigg|}dr\cr
&\leq& C_{T,|x_0|} \Lambda \bigg(1+\sup _{0 \leq s \leq t}\int_{0}^{s} e^{-\lambda (s-r)}\cr
&&\times\big[(s-r)^{-2\alpha} \|X^\varepsilon\|_{\lambda,t} +(s-r)^{-\alpha} \|X^\varepsilon\|_{1,\lambda,t}\big]dr\bigg)\cr
&\leq& K \Lambda(1+\lambda^{2 \alpha-1}\|X^{\varepsilon}\|_{\lambda, t}+\lambda^{\alpha-1}\|X^{\varepsilon}\|_{1, \lambda, t}),
\end{eqnarray}
where we have used the estimate
\begin{eqnarray}\label{inq2}
\int_{0}^{t} e^{-\lambda(t-r)}(t-r)^{-2 \alpha}d r &=&  \frac{1}{\lambda} \int_{0}^{\lambda t} e^{-q} \lambda^{2 \alpha} q^{-2 \alpha} d q \cr
&\leq& \lambda^{2 \alpha-1} \int_{0}^{\infty} e^{-q}q^{-2 \alpha} d q\cr
&\leq& C \lambda^{2 \alpha-1}.
\end{eqnarray}

Putting $\lambda=(4 K \Lambda)^{\frac{1}{1-\alpha}}$, we get from the inequality (\ref{in1}) that
\begin{eqnarray}\label{in4}
\|X^\varepsilon\|_{\lambda, t} \leq \frac{4}{3} K \Lambda(1+\lambda^{-1}\|X^\varepsilon\|_{1, \lambda, t}).
\end{eqnarray}
Then,
plugging this to the inequality (\ref{in2}) and making simple transformations, we arrive at
\begin{eqnarray*}
	\|X^\varepsilon\|_{1, \lambda, t} \leq \frac{3}{2} K \Lambda+2 (K \Lambda)^{1 /(1-\alpha)}\leq C_{T,|x_0|}\Lambda^{1 /(1-\alpha)}.
\end{eqnarray*}
Substituting this into (\ref{in4}), we get
$$
\|X^\varepsilon\|_{\lambda, t} \leq C_{T,|x_0|} \Lambda^{1 /(1-\alpha)}.
$$
Thus, we have
\begin{eqnarray*}
	\|X^\varepsilon\|_{\alpha, \infty} &\leq& e^{\lambda T}(\|X^\varepsilon\|_{\lambda, T}+\|X^\varepsilon\|_{1,\lambda, T}) \cr
	&\leq& C_{T,|x_0|} \exp(C_{T,|x_0|} \Lambda^{1 /(1-\alpha)})\Lambda^{1 /(1-\alpha)}\cr
	&\leq& C_{T,|x_0|} \exp(C_{T,|x_0|} \Lambda^{\frac{1}{1-\alpha}}_{\alpha}(B^{H}))(1+\Lambda^{\frac{1}{1-\alpha}}_{\alpha}(B^{H})).
\end{eqnarray*}
Since $0<\frac{1}{1-\alpha}<2$, by (\ref{expfbm}), we have $\mathbb{E}[\exp(\Lambda_{\alpha}^{\frac{1}{1-\alpha}}(B^H))]<\infty.$\\
This completed the proof of Lemma \ref{lemb}. \qed

Moreover, using similar techniques in Lemma \ref{lemb}, we can prove
\begin{eqnarray}\label{hatxbarx}
\quad \quad \|\hat{X}^\varepsilon\|_{\alpha, \infty}+\|\bar{X}\|_{\alpha, \infty}  \leq C_{T,|x_0|} \exp(C_{T,|x_0|} \Lambda^{\frac{1}{1-\alpha}}_{\alpha}(B^{H}))(1+\Lambda^{\frac{1}{1-\alpha}}_{\alpha}(B^{H})).
\end{eqnarray}
and
\begin{eqnarray}\label{barbound}
\mathbb{E}[\mathbb\|\bar X\|_{\alpha,\infty}^{p}] +\mathbb{E} [\|\hat{X}^\varepsilon\|_{\alpha,\infty}^{p}] \leq C_{\alpha,T,|x_0|}.
\end{eqnarray}
Here, we omit the proof.

\begin{lem}\label{lemregu}
	Suppose that {\rm (H1)-(H3)} hold. Then, if $0 \leq t\leq t+h\leq T$, we have
	$$\mathbb{E}[|X^\varepsilon_{t+h}-X^\varepsilon_{t}|^2]\leq C_{\alpha,\beta,T,|x_0|} h^{2-2\alpha}.$$
\end{lem}
\para{Proof:} From (\ref{orinteg1}), by (H1)-(H3), we have
\begin{eqnarray*}
	\mathbb{E}[|X^\varepsilon_{t+h}-X^\varepsilon_{t}|^2]&\leq &\mathbb{E}\bigg[\bigg|\int_{t}^{t+h} b_1 (r,X^\varepsilon_{r}, Y^\varepsilon_{r}) d r\bigg|^2\bigg]+\mathbb{E}\bigg[\bigg|\int_{t}^{t+h} \sigma_1(r,X^\varepsilon_{r})d B^H_{r}\bigg|^2\bigg]\cr
	&=:&A_1^{h}+A_2^{h}.
\end{eqnarray*}

Firstly, for $A_1^{h}$, by (H3), it is easy to get
$A_1^{h} \leq Ch^2.$ For the second term, by Remark \ref{fbmito} and (H1), we firstly give the following estimate:
\begin{eqnarray}\label{t-s}
\bigg|\int_{s}^{t} \sigma_1(r,X^\varepsilon_{r}) d B^H_{r}\bigg|&\leq& \Lambda_{\alpha}(B^H) \bigg(\int_{s}^{t} \frac{|\sigma_1(r,X^\varepsilon_{r})|}{(r-s)^{\alpha}}dr\cr
&&+\int_{s}^{t} \int_{s}^{r} \frac{|\sigma_1(r,X^\varepsilon_{r})-\sigma_1(q,X^\varepsilon_{q})|}{(r-q)^{1+\alpha}}dqdr\bigg)\cr
&\leq& \Lambda_{\alpha}(B^H) (1+\|X^{\varepsilon}\|_{\alpha,\infty})\cr
&&\times\bigg(\int_{s}^{t}(r-s)^{-\alpha}dr+\int_{s}^{t}[(r-s)^{\beta-\alpha}+1]dr\bigg)\cr
&\leq& C_{\alpha,\beta,T} \Lambda_{\alpha}(B^H) (1+\|X^{\varepsilon}\|_{\alpha,\infty})(t-s)^{1-\alpha}.
\end{eqnarray}

Then, using similar techniques, we can prove that
\begin{eqnarray}\label{t-s2}
&&\bigg|\int_{s}^{t} \sigma_1(r,\hat{X}^\varepsilon_{r}) d B^H_{r}\bigg|+\bigg|\int_{s}^{t} \sigma_1(r,\bar{X}_{r}) d B^H_{r}\bigg| \cr
&&\quad \quad \quad \quad \quad \leq C_{\alpha,\beta,T} \Lambda_{\alpha}(B^H) (1+\|\hat{X}^{\varepsilon}\|_{\alpha,\infty}+\|\bar{X}\|_{\alpha,\infty})(t-s)^{1-\alpha}.
\end{eqnarray}

To proceed, by Remark \ref{fbmito} and Lemma \ref{lemb}, we have
\begin{eqnarray*}
	A_2^{h} &\leq& C_{\alpha,\beta,T} \mathbb{E}[(\Lambda_{\alpha}(B^H) (1+\|X^{\varepsilon}\|_{\alpha,\infty}))^2] h^{2-2\alpha} \cr
	&\leq& C_{\alpha,\beta,T,|x_0|} h^{2-2\alpha}.
\end{eqnarray*}
Thus, we deduce the desired estimate.
\qed
\begin{lem}	\label{ybound}
	Suppose that {\rm (H1) (H2) and (H4)} hold.  Then, we have
	$$
	\sup_{t\in [0,T]}\mathbb{E}[|Y^\varepsilon_t|^{2}] \leq C_{\alpha,T,|x_0|,|y_0|}.
	$$
\end{lem}
\para{Proof:} Using It\^{o} formula, we have
\begin{eqnarray*}
	\mathbb{E}[|Y_{t}^{\varepsilon}|^{2}]=|y_0|^{2}+\frac{2}{\varepsilon} \mathbb{E} \bigg[\int_{0}^{t}\left\langle b_2(X_{s}^{\varepsilon}, Y_{s}^{\varepsilon}), Y_{s}^{\varepsilon}\right\rangle d s\bigg]+\frac{1}{\varepsilon} \mathbb{E} \bigg[\int_{0}^{t}\left |\sigma_{2}(X_{s}^{\varepsilon}, Y_{s}^{\varepsilon})\right|^{2} d s\bigg],
\end{eqnarray*}
then by (H4), we have
\begin{eqnarray*}
	\frac{d}{d t} \mathbb{E}[|Y_{t}^{\varepsilon}|^{2}]&=&\frac{2}{\varepsilon} \mathbb{E}[\langle b_2 (X_{t}^{\varepsilon}, Y_{t}^{\varepsilon}), Y_{t}^{\varepsilon}\rangle]+\frac{1}{\varepsilon} \mathbb{E}[|\sigma_{2}(X_{t}^{\varepsilon}, Y_{t}^{\varepsilon})|^{2}] \cr
	&\leq&-\frac{\beta_2}{\varepsilon} \mathbb{E}[|Y_{t}^{\varepsilon}|^{2}]+\frac{C}{\varepsilon} \mathbb{E}[|X_{t}^{\varepsilon}|^{2}]+\frac{C}{\varepsilon}.
\end{eqnarray*}
Hence, by Gronwall's inequality \cite[pp. 584]{Givon2007}, and Lemma \ref{lemb}, we obtain
\begin{eqnarray*}
	\mathbb{E}[|Y_{t}^{\varepsilon}|^{2}]& \leq& |y_0|^{2} e^{-\frac{\beta_2}{\varepsilon} t}+\frac{C}{\varepsilon} \int_{0}^{t} e^{-\frac{\beta_2}{\varepsilon}(t-s)}(1+\mathbb{E}[|X_{s}^{\varepsilon}|^{2}]) d s \cr  & \leq& C_{\alpha,T,|x_0|}(1+|y_0|^{2}).
\end{eqnarray*}
This completed the proof of Lemma \ref{ybound}. \qed

\para{Step 2: The estimates for $|Y^\varepsilon_t-\hat{Y}^\varepsilon_t|$ and $ \|X_{t}^{\varepsilon}-\hat{X}_{t}^{\varepsilon}\|_{\alpha}$.}
\begin{lem}\label{yhat} 	
	Suppose that {\rm (H1), (H2) and (H4)} hold.  Then, we have
	$$
	\sup_{t\in [0,T]}\mathbb{E} [|Y^\varepsilon_t-\hat{Y}^\varepsilon_t|^{2}]\leq  C_{\alpha,\beta,T,|x_0|} \delta.
	$$
\end{lem}
\para{Proof:} Using It\^{o} formula again, we have
\begin{eqnarray*}
	\mathbb{E}[|Y_{t}^{\varepsilon}-\hat{Y}^\varepsilon_t|^{2}]&=&\frac{2}{\varepsilon} \mathbb{E} \bigg[\int_{0}^{t}\langle b_2(X_{s}^{\varepsilon}, Y_{s}^{\varepsilon})-b_2(X_{s( \delta)}^{\varepsilon}, \hat{Y}_{s}^{\varepsilon}) , Y_{s}^{\varepsilon}-\hat{Y}^\varepsilon_s\rangle d s\bigg]\cr
	&&+\frac{1}{\varepsilon} \mathbb{E}\bigg[ \int_{0}^{t}|\sigma_{2}(X_{s}^{\varepsilon}, Y_{s}^{\varepsilon})-\sigma_{2}(X_{s( \delta)}^{\varepsilon}, \hat{Y}_{s}^{\varepsilon})|^{2} d s\bigg]\cr
	&=&\frac{1}{\varepsilon} \mathbb{E} \bigg[\int_{0}^{t}(2 \langle b_2(X_{s}^{\varepsilon}, Y_{s}^{\varepsilon})-b_2(X_{s}^{\varepsilon}, \hat{Y}_{s}^{\varepsilon}), Y_{s}^{\varepsilon}- \hat{Y}_{s}^{\varepsilon} \rangle \cr
	&&+|\sigma_{2}(X_{s}^{\varepsilon}, Y_{s}^{\varepsilon})-\sigma_{2}(X_{s}^{\varepsilon}, \hat{Y}_{s}^{\varepsilon})|^{2})d s\bigg]\cr
	&&+\frac{2}{\varepsilon} \mathbb{E} \bigg[\int_{0}^{t}\langle b_2(X_{s}^{\varepsilon}, \hat{Y}_{s}^{\varepsilon})-b_2(X_{s( \delta)}^{\varepsilon}, \hat{Y}_{s}^{\varepsilon}) , Y_{s}^{\varepsilon}- \hat{Y}_{s}^{\varepsilon} \rangle d s\bigg]\cr
	&&+\frac{2}{\varepsilon} \mathbb{E}\bigg[\int_{0}^{t} \langle \sigma_{2}(X_{s}^{\varepsilon},  Y_{s}^{\varepsilon})-\sigma_{2}(X_{s}^{\varepsilon}, \hat{Y}_{s}^{\varepsilon}),\cr
	&&  \quad \sigma_{2}(X_{s}^{\varepsilon}, \hat{Y}_{s}^{\varepsilon})-\sigma_{2}(X_{s( \delta)}^{\varepsilon}, \hat{Y}_{s}^{\varepsilon})\rangle d s\bigg]\cr
	&&+\frac{1}{\varepsilon} \mathbb{E}\bigg[ \int_{0}^{t}|\sigma_{2}(X_{s}^{\varepsilon}, \hat{Y}_{s}^{\varepsilon})-\sigma_{2}(X_{s( \delta)}^{\varepsilon}, \hat{Y}_{s}^{\varepsilon})|^{2} d s\bigg].
\end{eqnarray*}
By (H4) and Young's inequality, we have
\begin{eqnarray*}
	\frac{d}{dt}\mathbb{E}[|Y_{t}^{\varepsilon}-\hat{Y}^\varepsilon_t|^{2}]&\leq & \frac{-\beta_1}{\varepsilon} \mathbb{E}[ |Y_{t}^{\varepsilon}- \hat{Y}_{t}^{\varepsilon}|^2]+\frac{C}{\varepsilon} \mathbb{E}[ |X_{t}^{\varepsilon}-X_{t( \delta)}^{\varepsilon}||Y_{t}^{\varepsilon}- \hat{Y}_{t}^{\varepsilon}|]\cr
	&&+\frac{C}{\varepsilon} \mathbb{E} [|X_{t}^{\varepsilon}-X_{t( \delta)}^{\varepsilon}|^2]\cr
	&\leq & \frac{-\beta_1}{2\varepsilon} \mathbb{E}[ |Y_{t}^{\varepsilon}- \hat{Y}_{t}^{\varepsilon}|^2]+\frac{C}{\varepsilon} \mathbb{E} [|X_{t}^{\varepsilon}-X_{t( \delta)}^{\varepsilon}|^2].
\end{eqnarray*}

Then by Lemma \ref{lemregu} and Gronwall's inequality \cite[pp. 584]{Givon2007}, we have
\begin{eqnarray*}
	\mathbb{E}[|Y_{t}^{\varepsilon}-\hat{Y}^\varepsilon_t|^{2}]\leq C_{\alpha,\beta,T,|x_0|} \frac{
		\delta}{\varepsilon}\int_{0}^{t}e^{-\frac{\beta_1(t-s)}{2\varepsilon}}ds \leq C_{\alpha,\beta,T,|x_0|} \delta.
\end{eqnarray*}
This completed the proof of Lemma \ref{yhat}. \qed

\begin{lem}\label{x-xhat}
	Suppose that {\rm (H1)-(H4)} hold. Then, we have
	\begin{eqnarray*}
		\mathbb{E} \big[\sup _{t\in [0, T]}\|X_{t}^{\varepsilon}-\hat{X}_{t}^{\varepsilon}\|^2_{\alpha}\big] \leq C_{\alpha,\beta,T,|x_0|,|y_0|} \delta.
	\end{eqnarray*}
\end{lem}
\para{Proof:}
In order to estimate $\mathbb{X}:= \mathbb{E} \big[\sup _{t\in [0, T]}\|X_{t}^{\varepsilon}-\hat{X}_{t}^{\varepsilon}\|^2_{\alpha}\big]$, note that $A:=\big\| \int_{0}^{t} f(s)d s \big \|_{\alpha}$ needs to be estimated in advance, here $f:[0,T] \rightarrow\mathbb{R}^d$ is a measurable function. Using Fubini's theorem, we have
\begin{eqnarray}\label{falpha}
A
&\leq& \Big| \int_{0}^{t} f(s)d s \Big| +\int_{0}^{t}(t-s)^{-1-\alpha}\int_{s}^{t} |f(r)|d rds\cr
&\leq& t^\alpha \int_{0}^{t} (t-r)^{-\alpha} |f(r)|d r +C_{\alpha}\int_{0}^{t}(t-r)^{-\alpha}|f(r)|dr\cr
&\leq&C_{\alpha,T}\int_{0}^{t}(t-r)^{-\alpha}|f(r)|dr.
\end{eqnarray}

From (\ref{mildx}) and (\ref{xhat}), by (\ref{falpha}), H\"older's inequality, (H2), Lemmas \ref{lemregu} and \ref{yhat}, we have
\begin{eqnarray*}
	\mathbb{X}
	&\leq& C \mathbb{E}\bigg[\sup _{t\in [0, T]}\bigg\| \int_{0}^{t} (b_1(s,X^\varepsilon_s,Y^\varepsilon_{s}) -b_1(s,X^\varepsilon_s,\hat{Y}^\varepsilon_{s}) )d s \bigg \|_{\alpha}^2\bigg]\cr
	&&+C \mathbb{E}\bigg[\sup _{t\in [0, T]}\bigg\| \int_{0}^{t} (b_1(s,X^\varepsilon_s,\hat{Y}^\varepsilon_{s}) -b_1(s,X^\varepsilon_{s(\delta)},\hat{Y}^\varepsilon_{s}))d s \bigg \|_{\alpha}^2\bigg]\cr
	&&+C \mathbb{E}\bigg[\sup _{t\in [0, T]}\bigg\| \int_{0}^{t} (b_1(s,X^\varepsilon_{s(\delta)},\hat{Y}^\varepsilon_{s})- b_1(s(\delta),X_{s(\delta)}^{\varepsilon}, \hat{Y}_{s}^{\varepsilon}) )d s \bigg \|_{\alpha}^2\bigg]\cr
	&\leq& C_{\alpha,T}\int_{0}^{T}\mathbb{E}[ |b_1(s,X^\varepsilon_s,Y^\varepsilon_{s}) -b_1(s,X^\varepsilon_s,\hat{Y}^\varepsilon_{s})|^{2}]d s\cr
	&&+C_{\alpha,T}\int_{0}^{T} \mathbb{E}[ |b_1(s,X^\varepsilon_s,\hat{Y}^\varepsilon_{s}) -b_1(s,X^\varepsilon_{s(\delta)},\hat{Y}^\varepsilon_{s})|^{2}]d s\cr
	&&+C_{\alpha,T} \int_{0}^{T} \mathbb{E} [ |b_1(s,X^\varepsilon_{s(\delta)},\hat{Y}^\varepsilon_{s})- b_1(s(\delta),X_{s(\delta)}^{\varepsilon}, \hat{Y}_{s}^{\varepsilon})|^{2}]d s\cr
	&\leq& C_{\alpha,T} \int_{0}^{T}\mathbb{E}[|X^\varepsilon_s-X_{s(\delta)}^{\varepsilon}|^{2}]+\mathbb{E}[|Y^\varepsilon_{s}  - \hat{Y}_{s}^{\varepsilon}|^{2}]ds+C_{\alpha,T} \delta^2 \cr
	&\leq& C_{\alpha,\beta,T,|x_0|,|y_0|} \delta.
\end{eqnarray*}
This completed the proof of Lemma \ref{x-xhat}. \qed

\para{Step 3: The estimate for $ \|\bar{X}_{t}-\hat{X}_{t}^{\varepsilon} \|_{\alpha}$.}
\begin{lem}\label{x-xbar}
	Suppose that {\rm (H1)-(H4)} hold. Then, we have
	\begin{eqnarray*}
		\mathcal{A}&:=&\mathbb{E}\bigg[ \sup_{t\in [0, T]}\bigg\|\int_{0}^{t}(b_1(s(\delta), X_{s(\delta)}^{\varepsilon}, \hat{Y}_{s}^{\varepsilon})-\overline{b}_1(s(\delta), X_{s(\delta)}^{\varepsilon})) d s\bigg\|_{\alpha}^{2}\bigg] \cr
		&\leq& C_{\alpha,\beta,T,|x_0|,|y_0|} (\varepsilon \delta^{-1}+\delta).
	\end{eqnarray*}
\end{lem}
\para{Proof:} By elementary inequality,  we have
\begin{eqnarray*}
	\mathcal{A}&\leq & C \mathbb{E}\bigg[ \sup_{t\in [0, T]}\bigg|\sum_{k=0}^{\lfloor \frac{t}{\delta}\rfloor-1} \int_{k\delta}^{(k+1)\delta}(b_1(k\delta, X_{k\delta}^{\varepsilon}, \hat{Y}_{s}^{\varepsilon})-\overline{b}_1(k\delta, X_{k\delta}^{\varepsilon})) ds\bigg|^{2}\bigg]\cr
	&&+C \mathbb{E}\bigg[ \sup_{t\in [0, T]}\bigg|\int_{\lfloor \frac{t}{\delta}\rfloor \delta}^{t}(b_1(s(\delta), X_{s(\delta)}^{\varepsilon}, \hat{Y}_{s}^{\varepsilon})-\overline{b}_1(s(\delta), X_{s(\delta)}^{\varepsilon})) ds \bigg|^{2}\bigg]\cr
	&&+ C \mathbb{E}\bigg[ \sup_{t\in [0, T]}\bigg(\int_{0}^{t}\frac{|\int_{s}^{t}(b_1(r(\delta), X_{r(\delta)}^{\varepsilon}, \hat{Y}_{r}^{\varepsilon})-\overline{b}_1(r(\delta), X_{r(\delta)}^{\varepsilon})) d r|}{(t-s)^{1+\alpha}} ds\bigg)^2\bigg]\cr
	&=:& \sum_{i=1}^{3}\mathcal{A}_i.
\end{eqnarray*}	
For $\mathcal{A}_1$ and $\mathcal{A}_2$, by (H3), we have
\begin{eqnarray*}
	\sum_{i=1}^{2}\mathcal{A}_i
	&\leq & C \mathbb{E}\bigg[\sup_{t\in [0, T]}\lfloor \frac{t}{\delta}\rfloor \sum_{k=0}^{\lfloor \frac{t}{\delta}\rfloor-1} \bigg|\int_{k\delta}^{(k+1)\delta}(b_1(k\delta, X_{k\delta}^{\varepsilon}, \hat{Y}_{s}^{\varepsilon})-\overline{b}_1(k\delta, X_{k\delta}^{\varepsilon})) ds\bigg|^{2}\bigg]\cr
	&&+C_T \delta^2 \cr
	&\leq & \frac{C_T}{\delta^2}  \max_{0\leq k \leq \lfloor \frac{T}{\delta}\rfloor-1} \mathbb{E}\bigg[\bigg|\int_{k\delta}^{(k+1)\delta}(b_1(k\delta, X_{k\delta}^{\varepsilon}, \hat{Y}_{s}^{\varepsilon})-\overline{b}_1(k\delta, X_{k\delta}^{\varepsilon})) ds\bigg|^{2}\bigg]\cr
	&&+C_T\delta^2\cr
	&\leq &
	C_T \frac{\varepsilon^2}{\delta^2} \max_{0 \leq k \leq \lfloor \frac{T}{\delta}\rfloor-1}\int_{0}^{\frac{\delta}{\varepsilon}}\int_{\zeta}^{\frac{\delta}{\varepsilon}}\mathcal{J}_{k}(s,\zeta)dsd\zeta +C_T\delta^2,
\end{eqnarray*}	
where $0\leq \zeta \leq s\leq \frac{\delta}{\varepsilon}$, and
\begin{eqnarray}\label{jk1}
\mathcal{J}_{k}(s,\zeta) &=&\mathbb{E}[\langle b_1(k\delta,X^{\varepsilon}_{k\delta},\hat{Y}^{\varepsilon}_{s\varepsilon+k\delta})
-\bar{b}_1(k\delta, X^{\varepsilon}_{k\delta}),\cr
&&\quad b_1(k\delta,X^{\varepsilon}_{k\delta},\hat{Y}^{\varepsilon}_{\zeta \varepsilon+k\delta})
-\bar{b}_1(k\delta,X^{\varepsilon}_{k\delta})\rangle].
\end{eqnarray}

Then, for $\mathcal{A}_3$, by H\"{o}lder's inequality and the fact that $\alpha<\frac{1}{2}$, we have
\begin{eqnarray*}
	\mathcal{A}_3 & \leq & C \mathbb{E}\bigg[ \sup_{t\in [0, T]} \int_{0}^{t}\frac{ds}{(t-s)^{\frac{1}{2}+\alpha}}\cr
	&&\times  \int_{0}^{t}\frac{ |\int_{s}^{t}(b_1(r(\delta), X_{r(\delta)}^{\varepsilon}, \hat{Y}_{r}^{\varepsilon})-\overline{b}_1(r(\delta), X_{r(\delta)}^{\varepsilon})) d r|^{2}}{(t-s)^{\frac{3}{2}+\alpha}} ds\bigg]\cr
	&\leq &C_{\alpha,T}\mathbb{E}\bigg[ \sup_{t\in [0, T]} \int_{0}^{t}\frac{|\int_{s}^{t}(b_1(r(\delta), X_{r(\delta)}^{\varepsilon}, \hat{Y}_{r}^{\varepsilon})-\overline{b}_1(r(\delta), X_{r(\delta)}^{\varepsilon})) d r|^{2}}{(t-s)^{\frac{3}{2}+\alpha}} \mathbf{1}_{\ell^c} ds\bigg]\cr
	&&+C_{\alpha,T}\mathbb{E}\bigg[ \sup_{t\in [0, T]} \int_{0}^{t}\frac{|\int_{s}^{t}(b_1(r(\delta), X_{r(\delta)}^{\varepsilon}, \hat{Y}_{r}^{\varepsilon})-\overline{b}_1(r(\delta), X_{r(\delta)}^{\varepsilon})) d r|^{2}}{(t-s)^{\frac{3}{2}+\alpha}} \mathbf{1}_{\ell} ds\bigg]\cr
	&=:&\mathcal{A}_{31}+\mathcal{A}_{32},
\end{eqnarray*}	
where $\mathbf{1}_{\cdot}$ is an indicator function, $\ell:=\{t < (\lfloor \frac{s}{\delta}\rfloor+2)\delta\}$ and $\ell^c:=\{t \geq (\lfloor \frac{s}{\delta}\rfloor+2)\delta\}$.

By (H3) and the fact that
$\lfloor \lambda_1 \rfloor -\lfloor \lambda_2 \rfloor \leq \lambda_1-\lambda_2 +1, $ for $\lambda_1 \geq \lambda_2 \geq 0$, we have
\begin{eqnarray*}
	&&\mathcal{A}_{31}\cr
	&\leq &C_{\alpha,T} \mathbb{E}\bigg[ \sup_{t\in [0, T]} \int_{0}^{t}\frac{|\int_{s}^{(\lfloor \frac{s}{\delta}\rfloor+1)\delta}(b_1(r(\delta), X_{r(\delta)}^{\varepsilon}, \hat{Y}_{r}^{\varepsilon})-\overline{b}_1(r(\delta), X_{r(\delta)}^{\varepsilon})) dr|^{2}}{(t-s)^{\frac{3}{2}+\alpha}} \mathbf{1}_{\ell^c}ds\bigg]\cr
	&&+C_{\alpha,T}\mathbb{E}\bigg[ \sup_{t\in [0, T]} \int_{0}^{t}\frac{|\int_{\lfloor \frac{t}{\delta}\rfloor\delta}^{t}(b_1(r(\delta), X_{r(\delta)}^{\varepsilon}, \hat{Y}_{r}^{\varepsilon})-\overline{b}_1(r(\delta), X_{r(\delta)}^{\varepsilon})) dr|^{2}}{(t-s)^{\frac{3}{2}+\alpha}} \mathbf{1}_{\ell^c}ds\bigg]\cr
	&&+ C_{\alpha,T} \mathbb{E}\bigg[ \sup_{t\in [0, T]}\int_{0}^{t} \frac{(\lfloor \frac{t}{\delta}\rfloor-\lfloor \frac{s}{\delta}\rfloor-1)}{(t-s)^{\frac{3}{2}+\alpha}}\cr
	&&\times \sum_{k=\lfloor \frac{s}{\delta}\rfloor +1}^{\lfloor \frac{t}{\delta}\rfloor-1} \bigg|\int_{k\delta}^{(k+1)\delta }(b_1(k\delta, X_{k\delta}^{\varepsilon}, \hat{Y}_{r}^{\varepsilon})-\overline{b}_1(k\delta, X_{k\delta}^{\varepsilon})) dr\bigg|^{2}\mathbf{1}_{\ell^c}ds\bigg]\cr
	&\leq &C_{\alpha,T}  \sup_{t\in [0, T]}\bigg(\int_{0}^{t}(t-s)^{-\frac{1}{2}-\alpha} ((\lfloor \frac{s}{\delta}\rfloor+1)\delta-s)\mathbf{1}_{\ell^c}  ds\bigg)\cr
	&&+ C_{\alpha,T}  \sup_{t\in [0, T]}\bigg(\int_{0}^{t}(t-s)^{-\frac{1}{2}-\alpha} (t- \lfloor \frac{t}{\delta}\rfloor \delta)\mathbf{1}_{\ell^c}  ds\bigg)\cr
	&&+C_{\alpha,T}\delta^{-1}  \mathbb{E}\bigg[ \sup_{t\in [0, T]} \int_{0}^{t}(t-s)^{\frac{1}{2}-\alpha}\cr
	&&\quad \times \sum_{k=\lfloor \frac{s}{\delta}\rfloor +1}^{\lfloor \frac{t}{\delta}\rfloor-1}\bigg|\int_{k\delta}^{(k+1)\delta }(b_1(k\delta, X_{k\delta}^{\varepsilon}, \hat{Y}_{r}^{\varepsilon}) -\overline{b}_1(k\delta, X_{k\delta}^{\varepsilon})) dr\bigg|^{2} \mathbf{1}_{\ell^c}ds\bigg]\cr
	&\leq & C_{\alpha,T} \delta
	+ C_{\alpha,T} \frac{\varepsilon^2}{\delta^2} \max_{0 \leq k \leq\lfloor \frac{T}{\delta}\rfloor-1} \int_{0}^{\frac{\delta}{\varepsilon}}\int_{\zeta}^{\frac{\delta}{\varepsilon} }\mathcal{J}_{k}(s,\zeta)dsd\zeta.
\end{eqnarray*}

For $\mathcal{A}_{32}$, set $\jmath:=\{\lfloor \frac{t}{\delta}\rfloor >1\}$ and $\jmath^c:=\{\lfloor \frac{t}{\delta}\rfloor \leq 1\}$ by (H3) and the fact that $t -s< \lfloor \frac{s}{\delta}\rfloor \delta -s + 2\delta \leq  2\delta $,  we have
\begin{eqnarray*}
	&&\mathcal{A}_{32} \cr
	&\leq&  C_{\alpha,T}\mathbb{E}\bigg[\sup_{t\in [0, T]}\int_{0}^{t(\delta)-\delta}\frac{ |\int_{s}^{t}(b_1(r(\delta), X_{r(\delta)}^{\varepsilon}, \hat{Y}_{r}^{\varepsilon})-\overline{b}_1(r(\delta), X_{r(\delta)}^{\varepsilon})) d r|^{2}}{(t-s)^{\frac{3}{2}+\alpha}}\mathbf{1}_{\jmath \bigcap \ell} ds\bigg]\cr
	&&+ C_{\alpha,T}\mathbb{E}\bigg[ \sup_{t\in [0, T]}  \int_{t(\delta)-\delta}^{t}\frac{|\int_{s}^{t}(b_1(r(\delta), X_{r(\delta)}^{\varepsilon}, \hat{Y}_{r}^{\varepsilon})-\overline{b}_1(r(\delta), X_{r(\delta)}^{\varepsilon})) d r|^{2}}{(t-s)^{\frac{3}{2}+\alpha}}\mathbf{1}_{\jmath \bigcap \ell}ds\bigg]\cr
	&& + C_{\alpha,T}\mathbb{E}\bigg[ \sup_{t\in [0, T]}  \int_{0}^{t}\frac{ |\int_{s}^{t}(b_1(r(\delta), X_{r(\delta)}^{\varepsilon}, \hat{Y}_{r}^{\varepsilon})-\overline{b}_1(r(\delta), X_{r(\delta)}^{\varepsilon})) d r|^{2}}{(t-s)^{\frac{3}{2}+\alpha}}\mathbf{1}_{\jmath^c \bigcap \ell}ds\bigg]\cr
	&\leq&  C_{\alpha,T} \delta^2 \sup_{t\in [0, T]}\bigg( \int_{0}^{t(\delta)-\delta}(t-s)^{-\frac{3}{2}-\alpha}\mathbf{1}_{\jmath \bigcap \ell} ds\bigg)\cr
	&&+C_{\alpha,T} \sup_{t\in [0, T]}\bigg( \int_{t(\delta)-\delta}^{t}(t-s)^{\frac{1}{2}-\alpha}\mathbf{1}_{\jmath \bigcap \ell}ds\bigg)\cr
	&&+ C_{\alpha,T} \sup_{t\in [0, T]}\bigg( \int_{0}^{t}(t-s)^{\frac{1}{2}-\alpha}\mathbf{1}_{\jmath^c \bigcap \ell}ds\bigg)\cr
	&\leq&  C_{\alpha,T} \delta^{\frac{3}{2}-\alpha}.
\end{eqnarray*}

Thus, we have
\begin{eqnarray*}
	\mathcal{A} \leq C_{\alpha,T} \frac{\varepsilon^2}{\delta^2} \max_{0 \leq k \leq \lfloor \frac{T}{\delta}\rfloor-1}\int_{0}^{\frac{\delta}{\varepsilon}}\int_{\zeta}^{\frac{\delta}{\varepsilon}}\mathcal{J}_{k}(s,\zeta)dsd\zeta + C_{\alpha,T} \delta.
\end{eqnarray*}	

Now, by the construction of  $\hat{Y}^\varepsilon$ and a time shift transformation, for any fixed $k$ and $s\in[0,\delta]$, we have
\begin{eqnarray*}
	\hat{Y}_{s+k\delta}^\varepsilon &=& \hat Y_{k\delta}^\varepsilon+\frac{1}{\varepsilon}\int_{k\delta}^{k\delta+s}b_2(X_{k\delta}^\varepsilon,\hat{Y}_r^\varepsilon)dr+\frac{1}{\sqrt{\varepsilon}}\int_{k\delta}^{k\delta+s}\sigma_{2}(X_{k\delta}^\varepsilon,\hat{Y}_r^\varepsilon)dW_r\cr
	&=&\hat Y_{k\delta}^\varepsilon+\frac{1}{\varepsilon}\int_0^sb_2\big(X_{k\delta}^\varepsilon,\hat{Y}_{r+k\delta}^\varepsilon\big)dr+\frac{1}{\sqrt{\varepsilon}}\int_0^s\sigma_{2}\big(X_{k\delta}^\varepsilon,\hat{Y}_{r+k\delta}^\varepsilon\big)dW_r^{*},
\end{eqnarray*}
where $W_t^{*}=W_{t+k\delta}-W_{k\delta}$ is the shift version of $W_t$, and hence they have the same distribution. Let $\bar{W}$ be a Wiener process and independent of $W$. Construct a process $Y^{X_{k\delta}^\varepsilon,\hat Y_{k\delta}^\varepsilon}$ by means of
\begin{eqnarray}\label{yxy}
Y_{s/\varepsilon}^{X_{k\delta}^\varepsilon,\hat Y_{k\delta}^\varepsilon}&=&\hat Y_{k\delta}^\varepsilon+\int_0^{s/\varepsilon}b_2\big(X_{k\delta}^\varepsilon,Y_{r}^{X_{k\delta}^\varepsilon,\hat Y_{k\delta}^\varepsilon}\big)dr+\int_0^{s/\varepsilon}\sigma_{2}\big(X_{k\delta}^\varepsilon,Y_{r}^{X_{k\delta}^\varepsilon,\hat Y_{k\delta}^\varepsilon}\big)d\bar{W}_r\cr
&=&\hat Y_{k\delta}^\varepsilon+\frac{1}{\varepsilon}\int_0^{s}b_2\big(X_{k\delta}^\varepsilon,Y_{r/\varepsilon}^{X_{k\delta}^\varepsilon,\hat Y_{k\delta}^\varepsilon}\big)dr\cr
&&+\frac{1}{\sqrt{\varepsilon}}\int_0^{s}\sigma_{2}\big(X_{k\delta}^\varepsilon,Y_{r/\varepsilon}^{X_{k\delta}^\varepsilon, \hat Y_{k\delta}^\varepsilon}\big)d\bar{\bar{W}}^\varepsilon_r,
\end{eqnarray}
where $\bar{\bar{W}}^\varepsilon_t=\sqrt{\varepsilon}\bar{W}_{t/\varepsilon}$ is the scaled version of $\bar{W}_t$. Because both $W^{*}$ and $\bar{\bar{W}}$ are independent of $(X_{k\delta}^\varepsilon,\hat Y_{k\delta}^\varepsilon)$, by comparison, yields
\begin{eqnarray}\label{dist}
(X_{k\delta}^\varepsilon,\{\hat{Y}_{s+k\delta}^\varepsilon\}_{s\in[0,\delta)}\big)\sim\big(X_{k\delta}^\varepsilon,\{Y_{s/\varepsilon}^{X_{k\delta}^\varepsilon,\hat{Y}_{k\delta}^\varepsilon}\}_{s\in[0,\delta)}),
\end{eqnarray}
where $\sim$ denotes coincidence in distribution sense.

Thus, for  $s\in [0,\delta),$ from (\ref{jk1}), we have
\begin{eqnarray*}
	&&\mathcal{J}_{k}(s,\zeta)\cr
	&=&\mathbb{E}[\langle b_1(k\delta,X^{\varepsilon}_{k\delta},Y^{X^{\varepsilon}_{k\delta},\hat{Y}^{\varepsilon}_{k\delta}}_{s})
	-\bar{b}_1(k\delta,X^{\varepsilon}_{k\delta}), b_1(k\delta,X^{\varepsilon}_{k\delta},Y^{X^{\varepsilon}_{k\delta},\hat{Y}^{\varepsilon}_{k\delta}}_\zeta)
	-\bar{b}_1(k\delta,X^{\varepsilon}_{k\delta})\rangle].
\end{eqnarray*}

Now, we present a claim which will be proved in Appendix B.
\para{Claim:}
$$\mathcal{J}_{k}(s,\zeta)\leq C e^{-\frac{\beta_1}{2}(s-\zeta)}\mathbb{E}[(1+|X^{\varepsilon}_{k\delta}|^{2}+|\hat{Y}^{\varepsilon}_{k\delta}|^{2})]\leq C_{\alpha,\beta,T,|x_0|,|y_0|}e^{-\frac{\beta_1}{2}(s-\zeta)},$$
where $\beta_1$ is defined in (H4). Here, Lemmas \ref{lemb}, \ref{ybound} and \ref{yhat} were used for the last inequality.

Therefore, by choosing $\delta=\delta(\varepsilon)$ such that $\frac{\delta}{\varepsilon}$ is sufficiently large, we have
\begin{eqnarray*}
	\mathcal{A}&\leq & C_{\alpha,T} \frac{\varepsilon^2}{\delta^2} \max_{0 \leq k \leq \lfloor \frac{T}{\delta}\rfloor-1}\int_{0}^{\frac{\delta}{\varepsilon}}\int_{\zeta}^{\frac{\delta}{\varepsilon}}e^{-\frac{\beta_1}{2}(s-\zeta)}dsd\zeta + C_{\alpha,T} \delta\cr
	&\leq & C_{\alpha,\beta,T,|x_0|,|y_0|} \frac{\varepsilon^2}{\delta^2} (\frac{2}{\beta_1}\frac{\delta}{\varepsilon}-\frac{4}{\beta_1^{2}}+e^{\frac{-\beta_1}{2}\frac{\delta}{\varepsilon}}) + C_{\alpha,T} \delta\cr
	&\leq &
	C_{\alpha,\beta,T,|x_0|,|y_0|} (\varepsilon \delta^{-1}+\delta).
\end{eqnarray*}	
This completed the proof of Lemma \ref{x-xbar}.\qed

Then, for each $R>1$, we define the following stopping time $\tau_R$,
\begin{eqnarray}\label{stoptime}
\tau_R :=\inf \{t\geq 0:\|B^{H}\|_{1-\alpha,\infty, t}\geq R \} \wedge T.
\end{eqnarray}
\begin{lem}\label{ptau}
	The following inequality holds (see also {\rm \cite[Lemma 4.4]{Mishura2011}}):
	$$\mathbb{P}\big(\tau_R < T\big) \leq R^{-1}\mathbb{E}[\|B^{H}\|^2_{1-\alpha,\infty, T}],$$
	and $R^{-1}\mathbb{E}[\|B^{H}\|^2_{1-\alpha,\infty, T}]$ tends to $0$ when $R \rightarrow \infty$.
\end{lem}
\para{Proof:} By Chebyshev's inequality, we have
\begin{eqnarray*}
	\mathbb{P}\big(\tau_R < T\big) \leq \mathbb{P}\big(\|B^{H}\|_{1-\alpha,\infty, T} \geq R\big)
	\leq R^{-1}\mathbb{E}[\|B^{H}\|^2_{1-\alpha,\infty, T}].
\end{eqnarray*}
Because $\|B^{H}\|_{1-\alpha,\infty, T}$ has moments of all order, see Lemma 7.5 in Nualart and R\u{a}\c{s}canu \cite{Rascanu2002}, thus we have
$\lim_{R \rightarrow \infty} R^{-1}\mathbb{E}[\|B^{H}\|^2_{1-\alpha,\infty, T}]=0.$
\qed	
\begin{lem}\label{xbar-xhat}
	Suppose that {\rm (H1)-(H4)} hold.  Then, we have
	\begin{eqnarray*}
		\mathbb{E}\big[\sup_{t\in [0,T]}\|\hat{X}_{t}^{\varepsilon}-\bar{X}_{t}\|_{\alpha}^2\big] &\leq& C_{\alpha,T,|x_0|}\sqrt{R^{-1}\mathbb{E}[\|B^{H}\|^2_{1-\alpha,\infty, T}]}\cr
		&&+C_{\alpha,\beta,\gamma,T, R, |x_0|,|y_0|} ({\varepsilon \delta^{-1}}+{\delta}).
	\end{eqnarray*}
\end{lem}
\para{Proof:} From (\ref{xave}) and (\ref{xhat}), we have
\begin{eqnarray}\label{xhat-xbarstop}
\mathbb{E}\big[\sup_{t\in [0,T]}\|\hat{X}_{t}^{\varepsilon}-\bar{X}_{t}\|_{\alpha}^2\big] &\leq&  \mathbb{E}\big[\sup_{t\in [0,T]}\|\hat{X}_{t}^{\varepsilon}-\bar{X}_{t}\|_{\alpha}^2\mathbf{1}_{\{\tau_R < T\}}\big]\cr
&&+ \mathbb{E}\big[\sup_{t\in [0,T]}\|\hat{X}_{t}^{\varepsilon}-\bar{X}_{t}\|_{\alpha}^2\mathbf{1}_{\{\tau_R \geq T\}}\big].
\end{eqnarray}

For the first supremum in the right-hand side of inequality (\ref{xhat-xbarstop}), by H\"older's inequality, we have
\begin{eqnarray*}
\mathbb{E}\big[\sup_{t\in [0,T]}\|\hat{X}_{t}^{\varepsilon}-\bar{X}_{t}\|_{\alpha}^2\mathbf{1}_{\{\tau_R < T\}}\big] \leq \big(\mathbb{E}\big[\sup_{t\in [0,T]}\|\hat{X}_{t}^{\varepsilon}-\bar{X}_{t}\|_{\alpha}^4\big]\big)^{\frac{1}{2}}\mathbb{P}\big(\tau_R < T\big)^{\frac{1}{2}}.
\end{eqnarray*}
It follows from Lemma \ref{ptau} that $\mathbb{P}\big(\tau_R < T\big) \leq R^{-1}\mathbb{E}[\|B^{H}\|^2_{1-\alpha,\infty, T}]$. Then, by (\ref{barbound}), summing up all bounds we obtain
\begin{eqnarray}\label{j0}
\quad \quad \mathbb{E}\big[\sup_{t\in [0,T]}\|\hat{X}_{t}^{\varepsilon}-\bar{X}_{t}\|_{\alpha}^2\mathbf{1}_{\{\tau_R < T\}}\big]  \leq C_{\alpha,T,|x_0|}\sqrt{R^{-1}\mathbb{E}[\|B^{H}\|^2_{1-\alpha,\infty, T}]}.
\end{eqnarray}

Now, for $\lambda \geq 1$ and let
\begin{eqnarray*}
	\mathbf{A}:=\mathbb{E}\big[\sup_{t\in [0,T]}e^{-\lambda t} \|\hat{X}_{t}^{\varepsilon}-\bar{X}_{t}\|_{\alpha}^2\mathbf{1}_{D}\big],
\end{eqnarray*}
where $D :=\{\|B^{H}\|_{1-\alpha,\infty, T} \leq R\} $. Then, we return to the second supremum in the right-hand side of inequality (\ref{xhat-xbarstop}),
\begin{eqnarray*}
	\mathbf{A}
	&\leq&C \mathbb{E}\bigg[\sup_{t\in [0,T]} e^{-\lambda t}\bigg\|\int_{0}^{t} (b_1(s(\delta),X_{s(\delta)}^{\varepsilon}, \hat{Y}_{s}^{\varepsilon})-\bar{b}_1(s(\delta),X_{s(\delta)}^{\varepsilon}) ) d s\bigg\|_{\alpha}^2\mathbf{1}_{D} \bigg]\cr
	&&+C  \mathbb{E}\bigg[\sup_{t\in [0,T]} e^{-\lambda t}\bigg \|\int_{0}^{t} (\bar{b}_1(s(\delta),X_{s(\delta)}^{\varepsilon}) -\bar{b}_1(s,X_{s(\delta)}^{\varepsilon}) ) d s\bigg\|_{\alpha}^2\mathbf{1}_{D} \bigg]\cr
	&&+C \mathbb{E}\bigg[\sup_{t\in [0,T]} e^{-\lambda t}\bigg \|\int_{0}^{t} (\bar{b}_1(s, X_{s(\delta)}^{\varepsilon}) -\bar{b}_1(s,{X}_{s}^{\varepsilon}) ) d s\bigg\|_{\alpha}^2\mathbf{1}_{D} \bigg]\cr
	&&+C \mathbb{E}\bigg[\sup_{t\in [0,T]} e^{-\lambda t}\bigg \|\int_{0}^{t} (\bar{b}_1(s,{X}_{s}^{\varepsilon})-\bar{b}_1(s,\hat{X}^\varepsilon_{s})) d s\bigg\|_{\alpha}^2\mathbf{1}_{D} \bigg]\cr
	&&+C \mathbb{E}\bigg[\sup_{t\in [0,T]} e^{-\lambda t}\bigg \|\int_{0}^{t} (\bar{b}_1(s,\hat{X}^\varepsilon_{s})-\bar{b}_1(s,\bar{X}_{s}) ) d s\bigg\|_{\alpha}^2\mathbf{1}_{D} \bigg]\cr
	&&+C \mathbb{E}\bigg[\sup_{t\in [0,T]} e^{-\lambda t}\bigg \|\int_{0}^{t}( \sigma_{1}(s,\hat{X}^\varepsilon_{s})-\sigma_{1}(s,\bar{X}_{s}) )d B^{H}_{s}\bigg\|_{\alpha}^2\mathbf{1}_{D} \bigg]\cr
	&&+C \mathbb{E}\bigg[\sup_{t\in [0,T]} e^{-\lambda t}\bigg \|\int_{0}^{t}( \sigma_{1}(s,{X}_{s}^{\varepsilon})-\sigma_{1}(s,\hat{X}^\varepsilon_{s}) )d B^{H}_{s}\bigg\|_{\alpha}^2\mathbf{1}_{D} \bigg]\cr
	&=:&\sum_{i= 1}^{7}\mathbf{A}_i.
\end{eqnarray*}

By Lemma \ref{x-xbar}, we can estimate the term $\mathbf{A}_1$,
\begin{eqnarray}\label{j1}
\mathbf{A}_1
\leq C_{\alpha,\beta,T,|x_0|,|y_0|} ({\varepsilon \delta^{-1}}+{\delta}).
\end{eqnarray}

Then, by (H2), (\ref{falpha}), Lemma \ref{lemregu} and Lemma \ref{x-xhat}, it is easy to obtain
\begin{eqnarray}\label{j234}
\mathbf{A}_2+\mathbf{A}_3+\mathbf{A}_4
&\leq& C_{\alpha,T} \mathbb{E}\bigg[ \int_{0}^{T} |\bar{b}_1(s(\delta),X_{s(\delta)}^{\varepsilon}) -\bar{b}_1(s,X_{s(\delta)}^{\varepsilon})|^2d s \bigg]\cr
&&+C_{\alpha,T}  \mathbb{E}\bigg[\int_{0}^{T} |\bar{b}_1(s, X_{s(\delta)}^{\varepsilon}) -\bar{b}_1(s,{X}_{s}^{\varepsilon})|^2d s \bigg]\cr
&&+C_{\alpha,T}  \mathbb{E}\bigg[\int_{0}^{T} |\bar{b}_1(s,{X}_{s}^{\varepsilon})-\bar{b}_1(s,\hat{X}^\varepsilon_{s})|^2 d s \bigg]\cr
&\leq& C_{\alpha,\beta,T,|x_0|,|y_0|} {\delta}.
\end{eqnarray}
For $\mathbf{A}_5$, by (\ref{inq1}) and (\ref{falpha}), we have
\begin{eqnarray}\label{j5}
\mathbf{A}_5
&\leq& C_{\alpha,T}\mathbb{E}\bigg[\sup_{t\in [0,T]} e^{-\lambda t}\int_{0}^{t} (t-s)^{-2 \alpha}|\bar{b}_1(s,\hat{X}^\varepsilon_{s}) )  -\bar{b}_1(s,\bar{X}_{s})|^2\mathbf{1}_{D}d s\bigg]\cr
&\leq & C_{\alpha,T}  \mathbb{E}\bigg[\sup_{t\in [0,T]}\int_{0}^{t} e^{-\lambda (t-s)}(t-s)^{-2\alpha} e^{-\lambda s}|\hat{X}^\varepsilon_{s}-\bar{X}_{s}|^2\mathbf{1}_{D}d s\bigg]\cr
&\leq& C_{\alpha,T}\mathbb{E}\big[\sup_{t\in [0,T]} e^{-\lambda t}\|\hat{X}^\varepsilon_{t}-\bar{X}_{t}\|_{\alpha}^2\mathbf{1}_{D}\big] \sup_{t\in [0, T]} \int_{0}^{t} e^{-\lambda (t-r)}(t-r)^{-2 \alpha} d r\cr
&\leq& C_{\alpha,T} \lambda^{2\alpha-1} \mathbb{E}\big[\sup_{t\in [0,T]} e^{-\lambda t}\|\hat{X}^\varepsilon_{t}-\bar{X}_{t}\|_{\alpha}^2\mathbf{1}_{D}\big].
\end{eqnarray}

For $\mathbf{A}_6,\mathbf{A}_7$, we firstly give the basic estimate for $
\mathbf{B}:=\|\int_{0}^t f(r) d B^H_{r}\|_{\alpha}$ where $f:[0,T] \rightarrow\mathbb{R}^d$ is a measurable function. By (\ref{fbm}) and (\ref{fbmterm1}), it is easy to get
\begin{align}\label{fbmalpha}
\begin{split}
\mathbf{B}
\leq C_{\alpha, T} \Lambda_{\alpha}(B^H) \int_{0}^t ((t-r)^{-2\alpha}+r^{-\alpha}) \big(|f(r)|+\int_{0}^{r}\frac{|f(r)-f(q)|}{(r-q)^{1+\alpha}}dq\big) dr.
\end{split}
\end{align}

Next, by Lemma 7.1 in Nualart and R\u{a}\c{s}canu \cite{Rascanu2002}, we have
\begin{align}\label{inequ}
\begin{split}
|\sigma(t_1, x_1)-\sigma(t_2, x_2) &-\sigma(t_1, x_3)+\sigma(t_2, x_4)|
\\
\leq & C |x_1-x_2-x_3+x_4|+C |x_1-x_3||t_2-t_1|^{\beta}\\
&+C |x_1-x_3| (|x_1-x_2|^{\gamma}+|x_3-x_4|^{\gamma}).
\end{split}
\end{align}

Thus, by (\ref{inequ}), we have
\begin{eqnarray*}
	&&\mathbf{A}_6\cr
	&\leq& C_{\alpha,T,R} \mathbb{E}\bigg[\sup_{t\in [0,T]} \bigg|\int_{0}^{t} e^{-\lambda t}[(t-r)^{-2\alpha}+r^{-\alpha}]\|\sigma_{1}(r,\hat{X}^\varepsilon_{r})-\sigma_{1}(r,\bar{X}_{r})\|_{\alpha}\mathbf{1}_{D}d r\bigg|^2\bigg]\cr
	&\leq& C_{\alpha,T,R} \mathbb{E}\bigg[\sup_{t\in [0,T]} \bigg| \int_{0}^{t} e^{-\lambda t}[(t-r)^{-2\alpha}+r^{-\alpha}]\cr
	&&\times (1+\Delta(\hat{X}^\varepsilon_{r})+\Delta(\bar{X}_r))\|\hat{X}^\varepsilon_{r}-\bar{X}_{r}\|_{\alpha}\mathbf{1}_{D}  d r\bigg|^2\bigg],
\end{eqnarray*}
where $\Delta(\hat{X}^\varepsilon_{r})= \int_{0}^{r} \frac{|\hat{X}^\varepsilon_{r}-\hat{X}^\varepsilon_q|^\gamma}{(r-q)^{1+\alpha}}dq$ and $\Delta(\bar{X}_{r})= \int_{0}^{r} \frac{|\bar{X}_{r}-\bar{X}_q|^\gamma}{(r-q)^{1+\alpha}}dq$.

By (\ref{t-s2}), we have
\begin{eqnarray}\label{delta}
\Delta(\hat{X}^\varepsilon_{r})+\Delta(\bar{X}_r)
&\leq&C_{\alpha,\beta,T} \Lambda^\gamma_{\alpha}(B^H) (1+\|\hat{X}^{\varepsilon}\|_{\alpha,\infty})^\gamma \int_{0}^{t} (t-s)^{(1-\alpha)\gamma-1-\alpha}ds\cr
&&+C_{\alpha,\beta,T} \Lambda^\gamma_{\alpha}(B^H) (1+\|\bar{X}\|_{\alpha,\infty})^\gamma \int_{0}^{t} (t-s)^{(1-\alpha)\gamma-1-\alpha}ds\cr
&\leq&C_{\alpha,\beta,T} \Lambda^\gamma_{\alpha}(B^H) (1+\|X^{\varepsilon}\|^\gamma_{\alpha,\infty}+\|\bar{X}\|^\gamma_{\alpha,\infty})\frac{t^{(1-\alpha)\gamma-\alpha}}{(1-\alpha)\gamma-\alpha}\cr
&\leq&C_{\alpha,\beta,\gamma,T} \Lambda^\gamma_{\alpha}(B^H) (1+\|X^{\varepsilon}\|^\gamma_{\alpha,\infty}+\|\bar{X}\|^\gamma_{\alpha,\infty}).
\end{eqnarray}
Here, we use the fact that $(1-\alpha)\gamma-\alpha>0$, since $\alpha \in (0, \frac{\gamma}{2})$.

Then, by (\ref{hatxbarx}) and (\ref{delta}), under the condition that $\|B^{H}\|_{1-\alpha,\infty, T} \leq R$, there exists a constant $C_{\alpha,\beta,\gamma,T,R} $, such that
\begin{eqnarray}\label{delta2}
\Delta(\hat{X}^\varepsilon_{r})+\Delta(\bar{X}_r) \leq
C_{\alpha,\beta,\gamma,T,R,|x_0|}.
\end{eqnarray}

Thus, by (\ref{inq1}), (\ref{inq2}) and (\ref{delta2}), we obtain
\begin{eqnarray}\label{j6}
\mathbf{A}_6\leq C_{\alpha,\beta,\gamma,T,R,|x_0|,|y_0|}  \lambda^{2\alpha-1}\mathbb{E}\big[\sup_{t\in [0,T]} e^{-\lambda t}\|\hat{X}^\varepsilon_{t}-\bar{X}_{t}\|_{\alpha}^2\mathbf{1}_{D}\big].
\end{eqnarray}

Using the similar techniques, and by Lemma \ref{x-xhat}, we get
\begin{eqnarray}\label{j7}
\mathbf{A}_7 &\leq& C_{\alpha,\beta,\gamma,T,R,|x_0|,|y_0|} \lambda^{2\alpha-1}\mathbb{E}\big[\sup_{t\in [0,T]} e^{-\lambda t}\|{X}^\varepsilon_{t}-\hat{X}^\varepsilon_{t}\|_{\alpha}^2\mathbf{1}_{D}\big]\cr
&\leq & C_{\alpha,\beta,\gamma,T,R,|x_0|,|y_0|} {\delta}.
\end{eqnarray}

According to estimates (\ref{j1}), (\ref{j234}), (\ref{j5}), (\ref{j6}) and (\ref{j7}), we obtain that
\begin{eqnarray*}
	\mathbf{A}
	&\leq& C_{\alpha,\beta,\gamma,T,R,|x_0|,|y_0|}\lambda^{2\alpha-1}\mathbb{E}\big[\sup_{t\in [0,T]} e^{-\lambda t}\|\hat{X}^\varepsilon_{t}-\bar{X}_{t}\|_{\alpha}^2\mathbf{1}_{D}\big]\cr
	&&+C_{\alpha,\beta,T,|x_0|,|y_0|} ({\varepsilon \delta^{-1}}+{\delta})+C_{\alpha,\beta,T,|x_0|,|y_0|} {\delta}+
	C_{\alpha,\beta,\gamma,T,R,|x_0|,|y_0|} {\delta}.
\end{eqnarray*}
Taking $\lambda$ large enough, such that $C_{\alpha,\beta,\gamma,T,R,|x_0|,|y_0|}\lambda^{2\alpha-1}<1$, we have
\begin{eqnarray}\label{J8}
\mathbb{E}\big[\sup_{t\in [0,T]} e^{-\lambda t}\|\hat{X}^\varepsilon_{t}-\bar{X}_{t}\|_{\alpha}^2\mathbf{1}_{D}\big]
\leq
C_{\alpha,\beta,\gamma,T,R,|x_0|,|y_0|} ({\varepsilon \delta^{-1}}+{\delta}).
\end{eqnarray}

Finally, by (\ref{j0}) and (\ref{J8}),  we obtain that
\begin{eqnarray*}
	\mathbb{E}\big[\sup_{t\in [0,T]} e^{-\lambda t}\|\hat{X}^\varepsilon_{t}-\bar{X}_{t}\|_{\alpha}^2\big]& \leq& C_{\alpha,T,|x_0|}\sqrt{R^{-1}\mathbb{E}[\|B^{H}\|^2_{1-\alpha,\infty, T}]}\cr
	&&+C_{\alpha,\beta,\gamma,T,R,|x_0|,|y_0|} ({\varepsilon \delta^{-1}}+{\delta}).
\end{eqnarray*}
Then, the statement follows. \qed

\para{Step 4: The estimate for $ \|\bar{X}_{t}-{X}_{t}^{\varepsilon}\|_{\alpha}$.}
By Lemma \ref{x-xhat} and Lemma \ref{xbar-xhat}, we have
\begin{eqnarray*}
	\mathbb{E}\big[\sup_{t\in [0,T]}\|{X}_{t}^{\varepsilon}-\bar{X}_{t}\|_{\alpha} ^2\big]  &\leq&\mathbb{E}\big[\sup_{t\in [0,T]}\|{X}_{t}^{\varepsilon}-\hat{X}_{t}^{\varepsilon}\|_{\alpha}^2\big]+\mathbb{E}\big[\sup_{t\in [0,T]}\|\hat{X}_{t}^{\varepsilon}-\bar{X}_{t}\|_{\alpha}^2\big]  \cr
	&\leq& C_{\alpha,T,|x_0|}\sqrt{R^{-1}\mathbb{E}[\|B^{H}\|^2_{1-\alpha,\infty, T}]}\cr
	&&+C_{\alpha,\beta,\gamma,T,R,|x_0|,|y_0|} ({\varepsilon \delta^{-1}}+{\delta}).
\end{eqnarray*}

Thus, choose $\delta=\varepsilon \sqrt{-\ln \varepsilon}$, we obtain
$$\limsup\limits_{\varepsilon \rightarrow 0}\mathbb{E}\big[\sup_{t\in [0,T]}\|{X}_{t}^{\varepsilon}-\bar{X}_{t}\|_{\alpha} ^2\big] \leq C_{\alpha,T,|x_0|}\sqrt{R^{-1}\mathbb{E}[\|B^{H}\|^2_{1-\alpha,\infty, T}]}.$$
Then, let $R \rightarrow \infty$ and by Lemma \ref{ptau}, we have
\begin{eqnarray*}
	\lim\limits_{\varepsilon \rightarrow 0}\mathbb{E}\big[\sup_{t\in [0,T]}\|{X}_{t}^{\varepsilon}-\bar{X}_{t}\|_{\alpha} ^2\big] =0.
\end{eqnarray*}
This completed the proof of Theorem \ref{thm1}. \qed

\appendix
\renewcommand\thesection{\arabic{section}}
\section*{Appendix A: Ergodicity}
For fixed $x\in \mathbb{R}^{d_1}$, consider the problem associated to fast motion (\ref{frozon}) with frozen slow component. If (H2) holds, then it is easy to prove for any fixed $x\in \mathbb{R}^{d_1}$ and initial value $y\in \mathbb{R}^{d_2}$,   (\ref{frozon}) has a unique strong solution denoted by $\{Y_t^{x,y}\}_{t\geq 0}$, which is a Markov process.

Let $\{P_t^x\}_{t\geq 0}$ be the transition semigroup of $\{Y_t^{x,y}\}_{t\geq 0}$, i.e. for any bounded measurable function $\varphi:\mathbb{R}^{d_2}\rightarrow \mathbb{R},$
$$P_{s}^x\varphi(y):=\mathbb{E}[\varphi(Y^{x,y}_s)],\quad y\in \mathbb{R}^{d_2},s\geq 0.$$
Under the assumption (H4), it is easy to prove that
\begin{eqnarray}\label{a1}
\mathbb{E}[|Y^{x,y}_t|^{2}] \leq e^{-\beta_2 t}|y|^2+C(1+|x|^2),
\end{eqnarray}
and $\{P_t^x\}_{t\geq 0}$ has a unique invariant measure $\mu^x$ (see \cite[Lemma 3.6 and Proposition 3.8]{Liu2019} for example) satisfying
\begin{eqnarray}\label{a2}
\int_{\mathbb{R}^{d_2}}|y|^k \mu^{x}(dy)\leq C (1+|x|^k),
\end{eqnarray}
for some $k\geq 1$.

\begin{lem}\label{y1-y2}
	Suppose that {\rm (H2)} and {\rm (H4)} hold. For any given value $x\in\mathbb{R}^{d_1},y_1,y_2\in \mathbb{R}^{d_2}$ and $t\geq 0$, we have
	\begin{eqnarray*}
		\mathbb{E}[|Y_{t}^{x,y_1}-Y_{t}^{x,y_2}|^{2}]\leq C e^{-\beta_1t}|y_1-y_2|^2.
	\end{eqnarray*}
\end{lem}
\para{Proof:} Using It\^{o} formula again, we have
\begin{eqnarray*}
	\mathbb{E}[|Y_{t}^{x,y_1}-Y_{t}^{x,y_2}|^{2}]&=&|y_1-y_2|^2\cr
	&&+\mathbb{E} \bigg[\int_{0}^{t}2 \langle b_2(x, Y_{s}^{x,y_1})-b_2(x,Y_{s}^{x,y_2}),Y_{s}^{x,y_1}-Y_{s}^{x,y_2} \rangle d s\bigg]\cr
	&&+\mathbb{E} \bigg[\int_{0}^{t}|\sigma_{2}(x, Y_{s}^{x,y_1})-\sigma_{2}(x,Y_{s}^{x,y_2})|^2ds\bigg].
\end{eqnarray*}

By (H4) and Gronwall's inequality \cite[pp. 584]{Givon2007}, we obtain
\begin{eqnarray*}
	\mathbb{E}[|Y_{t}^{x,y_1}-Y_{t}^{x,y_2}|^{2}]\leq C e^{-\beta_1t}|y_1-y_2|^2.
\end{eqnarray*}
This completed the proof. \qed

The estimate (\ref{a1}) and the classical Bogoliubov-Krylov argument imply the existence of invariant measures. For the uniqueness, by the estimate (\ref{a2}) and Lemma \ref{y1-y2}, it is sufficient to prove the following Lemma \ref{invariant1} (see also \cite[Proposition 3.8 and Proposition 3.9]{Liu2019} for example).
\begin{lem}\label{invariant1}
	Suppose that {\rm (H2)} and {\rm (H4)} hold. For any given value $x\in\mathbb{R}^{d_1},y\in \mathbb{R}^{d_2}$, there exist $C>0$ and $\beta_1>0$ such that for any Lipschitz function $\varphi: \mathbb{R}^{d_2} \rightarrow \mathbb{R}$,
	$$
	\bigg|P_{s}^{x} \varphi(y)-\int_{\mathbb{R}^{d_2}} \varphi(z) \mu^{x}(d z)\bigg| \leq C(1+|x|+|y|) e^{-\beta_1 s}|\varphi|_{L i p}, \quad s \geq 0,
	$$
	where $|\varphi|_{L i p}=\sup _{x \neq y} \frac{|\varphi(x)-\varphi(y)|}{|x-y|}.$	
 Moreover,
	\begin{eqnarray*}
	\quad \quad \bigg|\mathbb{E}[b_1(t,x,Y^{x,y}_s)] -\int_{\mathbb{R}^{d_2}} b_1(t,x,z)\mu^{x}(dz)\bigg|^2\leq C e^{-\beta_1 s}(1+|x|^2+|y|^2).
	\end{eqnarray*}
\end{lem}	

\begin{lem}\label{x1-x2}
	Suppose that {\rm (H2)} and {\rm (H4)} hold. For any given value $x_1,x_2\in\mathbb{R}^{d_1},y\in \mathbb{R}^{d_2}$, we have
	\begin{eqnarray*}
		\mathbb{E}[|Y_{t}^{x_1,y}-Y_{t}^{x_2,y}|^{2}]\leq C |x_1-x_2|^2.
	\end{eqnarray*}
\end{lem}
\para{Proof:} Using It\^{o} formula again and Young's inequality, we have
\begin{eqnarray*}
	\frac{d}{dt}\mathbb{E}[|Y_{t}^{x_1,y}-Y_{t}^{x_2,y}|^{2}]&=&\mathbb{E} [2\left\langle b_2(x_1, Y_{t}^{x_1,y})-b_2(x_2,Y_{t}^{x_2,y}),Y_{t}^{x_1,y}-Y_{t}^{x_2,y} \right\rangle \cr
	&&+ |\sigma_{2}(x_1, Y_{t}^{x_1,y})-\sigma_{2}(x_2,Y_{t}^{x_2,y})|^2]\cr
	&=&\mathbb{E} [2\left\langle b_2(x_1, Y_{t}^{x_1,y})-b_2(x_1,Y_{t}^{x_2,y}),Y_{t}^{x_1,y}-Y_{t}^{x_2,y} \right \rangle \cr
	&&+ |\sigma_{2}(x_1, Y_{t}^{x_1,y})-\sigma_{2}(x_1,Y_{t}^{x_2,y})|^2]\cr
	&&+\mathbb{E} [2\left\langle b_2(x_1, Y_{t}^{x_2,y})-b_2(x_2,Y_{t}^{x_2,y}),Y_{t}^{x_1,y}-Y_{t}^{x_2,y} \right\rangle \cr
	&&+ |\sigma_{2}(x_1, Y_{t}^{x_2,y})-\sigma_{2}(x_2,Y_{t}^{x_2,y})|^2]\cr
	&&+ \mathbb{E} [2\langle \sigma_{2} (x_1, Y_{t}^{x_1,y})-\sigma_2(x_1,Y_{t}^{x_2,y}), \cr
	&& \quad \sigma_{2} (x_1, Y_{t}^{x_2,y})-\sigma_2(x_2,Y_{t}^{x_2,y}) \rangle ]\cr
	&\leq& -\frac{\beta_1}{2} |Y_{t}^{x_1,y}-Y_{t}^{x_2,y}|^2+C|x_1-x_2|^2.
\end{eqnarray*}

By (H4) and Gronwall's inequality \cite[pp. 584]{Givon2007}, we obtain
\begin{eqnarray*}
	\mathbb{E}[|Y_{t}^{x_1,y}-Y_{t}^{x_2,y}|^{2}]\leq C |x_1-x_2|^2.
\end{eqnarray*}
This completed the proof. \qed

\section*{Appendix B: The Proof of Claim in Lemma \ref{x-xbar}}
Let $\bar W$ be as in (\ref{yxy}) and $\mathbb{Q}^{y}$ denote the probability law of the diffusion process $\{Y^{x}_{t}\}_{t\geq 0}$ which is governed by following equation
\begin{eqnarray*}
	dY^{x}_t=b_2(x,Y^{x}_t)dt+\sigma_{2}(x,Y^{x}_t)d\bar{W}_t,
\end{eqnarray*}
with initial value $Y^x_0=y$ and we denote the solution by $\{Y_t^{x,y}\}_{t\geq 0}$.
The expectation with respect to $\mathbb{Q}^{y}$ is denoted by $\mathbb{E}^{y}$. Hence, we have
$\mathbb{E}^{y}[\Psi(Y^{x}_t)]=\mathbb{E}[\Psi(Y^{x,y}_t)],$
for all bounded function $\Psi$. For more details on $\mathbb{Q}^{y}$, the readers are referred to \cite[pp.110]{oksendal2003stochastic}. Let $\mathscr{F}_{t}^{x}$ be the $\sigma$-field generated by $\{Y^{x,y} _{r},r\leq t\}$
and set
$$\mathcal{J}_{k}(s,\zeta,x,y)=\mathbb{E}[\langle b_1(k\delta, x,Y^{x,y}_s)-\bar{b}_1(k\delta, x), b_1(k\delta, x,Y^{x,y}_\zeta)-\bar{b}_1(k\delta, x)\rangle].$$
Then, we have
\begin{eqnarray*}
	\mathcal{J}_{k}(s,\zeta,x,y)
	&=&\mathbb{E}^{y}[\langle b_1(k\delta, x,Y^{x}_s)-\bar{b}_1(k\delta,x), b_1(k\delta,x,Y^{x}_\zeta)-\bar{b}_1(k\delta,x)\rangle]\cr
	&=&\mathbb{E}^{y}[\mathbb{E}^{y}[\langle b_1(k\delta,x,Y^{x}_s)-\bar{b}_1(k\delta,x), (b_1(k\delta,x,Y^{x}_\zeta)-\bar{b}_1(k\delta,x)\rangle|\mathscr{F}_{\zeta}^{x}]]\cr
	&=&\mathbb{E}^{y}[\langle b_1(k\delta,x,Y^{x}_\zeta)-\bar{b}_1(k\delta,x), \mathbb{E}^{y}[ (b_1(k\delta,x,Y^{x}_s)-\bar{b}_1(k\delta,x))|\mathscr{F}_{\zeta}^{x}]\rangle].
\end{eqnarray*}
To proceed, by invoking the Markov property of $\{Y^{x,y}_t\}_{t\geq 0}$, we have
\begin{eqnarray*}
	\mathcal{J}_{k}(s,\zeta,x,y)=
	\mathbb{E}^{y}[\langle b_1(k\delta,x,Y^{x}_\zeta)-\bar{b}_1(k\delta,x),  \mathbb{E}^{Y^{x,y}_\zeta}[b_1(k\delta,x,Y^{x}_{s-\zeta})-\bar{b}_1(k\delta,x)]\rangle ],	
\end{eqnarray*}
where $\mathbb{E}^{Y^{x,y}_\zeta}[b_1(k\delta,x,Y^{x}_{s-\zeta})-\bar{b}_1(k\delta,x)]$ means the function $$\mathbb{E}^{y}[b_1(k\delta,x,Y^{x}_{s-\zeta})-\bar{b}_1(k\delta,x)]$$ evaluated at $y=Y^{x,y}_\zeta$.

Using H\"{o}lder's inequality and the boundedness of the function $b_1$, we obtain
\begin{eqnarray*}
	\mathcal{J}_{k}(s,\zeta,x,y) &\leq& C(
	\mathbb{E}^{y}[|b_1(k\delta,x,Y^{x}_\zeta)-\bar{b}_1(k\delta,x)|^{2}])^{\frac{1}{2}}\cr
	&&\times(
	\mathbb{E}^{y}[|\mathbb{E}^{Y^{x,y}_\zeta}[b_1(k\delta,x,Y^{x}_{s-\zeta})-\bar{b}_1(k\delta, x)]|^{2}])^{\frac{1}{2}}.
\end{eqnarray*}
In view of Lemma \ref{invariant1}, we have
\begin{eqnarray}\label{jk}
\mathcal{J}_{k}(s,\zeta,x,y)\leq C (1+|x|^{2}+|y|^{2})e^{-\frac{\beta_1}{2}(s-\zeta)}.
\end{eqnarray}
Let $\mathscr{M}_{k\delta}^{\varepsilon}$ be the $\sigma$-field generated by $X^{\varepsilon}_{k\delta}$ and $\hat{Y}^{\varepsilon}_{k\delta}$ that is independent of $\{Y^{x,y}_{r}\}_{r\geq 0}$. By adopting the approach in \cite[Theorem 7.1.2]{oksendal2003stochastic} . We can show
\begin{eqnarray*}
	\mathcal{J}_{k}(s,\zeta)&=&\mathbb{E}[\mathbb{E}[\langle b_1(k\delta,X^{\varepsilon}_{k\delta},Y^{X^{\varepsilon}_{k\delta},\hat{Y}^{\varepsilon}_{k\delta}}_s)-\bar{b}_1(k\delta,X^{\varepsilon}_{k\delta}),\cr
	&&\quad \quad b_1(k\delta,X^{\varepsilon}_{k\delta},Y^{X^{\varepsilon}_{k\delta},\hat{Y}^{\varepsilon}_{k\delta}}_\zeta)-\bar{b}_1(k\delta,X^{\varepsilon}_{k\delta})\rangle|\mathscr{M}_{k\delta}^{\varepsilon}]]\cr
	&=&\mathbb{E}[\mathcal{J}_{k}(s,\zeta,x,y)|_{(x,y)=(X^{\varepsilon}_{k\delta},\hat{Y}^{\varepsilon}_{k\delta})}],
\end{eqnarray*}
which, with the aid of (\ref{jk}), yields
\begin{eqnarray*}
	\mathcal{J}_{k}(s,\zeta)
	\leq C\mathbb{E}[
	(1+|X^{\varepsilon}_{k\delta}|^{2}+|\hat{Y}^{\varepsilon}_{k\delta}|^{2})]e^{-\frac{\beta_1}{2}(s-\zeta)}.
\end{eqnarray*}
This completes  the proof of the claim. \qed

\section*{Acknowledgement}
B. Pei was partially supported by National Natural Science Foundation (NSF) of China under Grants No.11802216 and No.12172285, NSF of Chongqing under Grant No.cstc2021jcyj-msxmX0296, Fundamental Research Funds for the Central Universities, Young Talent fund of University Association for Science and Technology in Shaanxi, China, and JSPS Grant-in-Aid for JSPS Fellows under Grant No.JP18F18314. Y. Inahama was partially supported by JSPS KAKENHI under Grant No.JP15K04922 and Grant-in-Aid for JSPS Fellows under Grant No.JP18F18314. Y. Xu was partially supported by NSF of China under Grant No.12072264, Key International (Regional) Joint Research
Program of NSF of China under Grant No.12120101002, Research Funds for Interdisciplinary Subject of Northwestern Polytechnical University, and Shaanxi Provincial Key R\&D Program under Grants No.2019TD-010 and No.2020KW-013. B. Pei would like to thank JSPS for Postdoctoral Fellowships for Research in Japan (Standard).


\end{document}